\documentclass[aos,MSNbibl,dvips]{arximspdf}
\usepackage{graphicx}

\doi{10.1214/12-AOS1046} 
\volume{40}
\issue{5}
\pubyear{2012}
\firstpage{2667}
\lastpage{2696}

\makeatletter
\def\sff{\mathsf{f}}
\def\sfh{\mathsf{h}}
\def\sfg{\mathsf{g}}
\def\sfl{\mathsf{l}}
\def\hatthetak{\widehat{\theta}_\bk}
\def\gmin{g_{\min}}
\def\bX{\mathbf{X}}
\def\supp{\mathop{\operatorname{supp}}}
\def\1{\mathbf{1}}
\def\RR{\mathbb{R}}
\def\PP{\mathbb{P}}
\def\NN{\mathbb{N}}
\def\ZZ{\mathbb{Z}}
\def\ba{\mathbf{a}}

\def\gammaL{\gamma_{\vartheta}}

\def\gammaLL{\gamma_{L}}

\def\oomega{\bolds{\omega}}
\def\E{\mathbb{E}}
\def\Pb{\mathbf{P}}
\def\Ex{\mathbf{E}}
\def\mm{\mathbf{m}}

\def\bY{\mathbf{Y}}
\def\bk{\mathbf{k}}
\def\KL{\mathcal{K}}
\def\llambda{\bolds{\lambda}}
\def\bQ{\mathbf{Q}}
\def\zb{z_\gamma}

\newtheorem{theorem}{Theorem}
\newtheorem{lemma}{Lemma}
\newtheorem{corollary}[lemma]{Corollary}
\newtheorem{proposition}[lemma]{Proposition}
\newproclaim{remark}{Remark}

\makeatother

\begin{document}
\begin{frontmatter}

\title{Tight conditions for consistency of variable selection in the
context of high dimensionality\thanksref{T1}}
\thankstext{T1}{Supported by the French Agence Nationale de la
Recherche (ANR) under the grant PARCIMONIE.}
\runtitle{Consistent variable selection in nonparametric regression}

\begin{aug}
\author[A]{\fnms{La\"etitia} \snm{Comminges}\ead[label=e1]{laetitia.comminges@imagine.enpc.fr}}
\and
\author[A]{\fnms{Arnak S.} \snm{Dalalyan}\corref{}\ead[label=e2]{dalalyan@imagine.enpc.fr}}
\runauthor{L. Comminges and A. S. Dalalyan}
\affiliation{Universit\'e Paris Est/ENPC and ENSAE-CREST}
\address[A]{LIGM/IMAGINE\\
Universit\'e Paris Est/ENPC and ENSAE-CREST\\
77455 Marne-la-Vall\'ee cedex 2\\
France\\
\printead{e1}\\
\phantom{E-mail:\ }\printead*{e2}} 
\end{aug}

\received{\smonth{6} \syear{2011}}
\revised{\smonth{3} \syear{2012}}

%
\begin{abstract}
We address the issue of variable selection in the regression model with
very high ambient dimension, that is,
when the number of variables is very large. The main focus is on the
situation where the number of relevant variables, called
intrinsic dimension, is much smaller than the ambient dimension~$d$.
Without assuming any parametric form of the underlying
regression function, we get tight conditions making it possible to
consistently estimate the set of relevant variables.
These conditions relate the intrinsic dimension to the ambient
dimension and to the sample size. The procedure that is provably
consistent under these tight conditions is based on comparing quadratic
functionals of the empirical Fourier coefficients with
appropriately chosen threshold values.

The asymptotic analysis reveals the presence of two quite different re\-
gimes. The first regime is when the intrinsic dimension is fixed.
In this case the situation in nonparametric regression is the same as
in linear regression, that is, consistent variable selection
is possible if and only if $\log d$ is small compared to the sample
size $n$. The picture is different in the second regime, that is,
when the number of relevant variables denoted by $s$ tends to infinity
as $n\to\infty$. Then we prove that consistent variable selection
in nonparametric set-up is possible only if $s+\log\log d$ is small
compared to $\log n$. We apply these results to derive minimax separation
rates for the problem of variable selection.
\end{abstract}

%


\end{frontmatter}

\section{Introduction}\label{sec1}

Real-world data such as those obtained from neuroscience, chemometrics,
data mining or sensor-rich
environments are often extremely high-dimensional, severely
underconstrained (few data samples compared
to the dimensionality of the data) and interspersed with a large number
of irrelevant or redundant
features. Furthermore, in most situations the data is contaminated by
noise, making it even more difficult
to retrieve useful information from the data. Relevant variable
selection is a compelling approach for
addressing statistical issues in the scenario of high-dimensional and
noisy data with small sample size.
Starting from Mallows~\cite{Mallows}, Akaike~\cite{Akaike}, Schwarz
\cite{Schwarz} who introduced,
respectively, the famous criteria $C_p$, AIC and
BIC, the problem of variable selection was extensively studied in the
statistical
and machine learning literature both from the theoretical and
algorithmic viewpoints. It appears, however, that
the theoretical limits of performing variable selection in the context
of nonparametric regression are still
poorly understood, especially when the number of variables, denoted by
$d$ and referred to as
ambient dimension, is much larger than the sample size $n$. The purpose
of the present work is to explore
this setting under the assumption that the number of relevant
variables, hereafter called intrinsic
dimension and denoted by $d^*$, may grow with the sample size but
remains much smaller than $d$.

In the important particular case of linear regression, the latter
scenario was the subject of a number of recent
studies. Many of them rely on $\ell_1$-norm penalization~\cite
{Tibsh,Zhao,Meinshausen} and
constitute an attractive alternative to iterative variable selection
procedures~\cite{alquier,Zhang09}
and to marginal regression or correlation screening \cite
{Wass09,Fan09}. Promising results for feature selection
are also obtained by conformal prediction~\cite{Hebiri2010},
(minimax) concave penalties~\cite{Fan01,Fan11,ZhangCH10},
Bayesian approach~\cite{Scott10} and higher criticism \cite
{Donoho09}. Extensions to other settings including logistic
regression, generalized linear model and Ising model were carried out
in~\cite{Bunea09,Ravik,Fan09}, respectively.
Variable selection in the context of groups of variables with disjoint
or overlapping groups was studied
by~\cite{Yuan06,svssin,Pontil,Obozinski,Huang10}. Hierarchical
procedures for selection of relevant variables were proposed
by~\cite{Bach09,Bickeletal,BinYu09}.

It is now well understood that in the Gaussian sequence model and in
the high-dimensional linear regression with a
Gram matrix satisfying some variant of irrepresentable condition,
consistent estimation of the pattern of relevant
variables---also called the sparsity pattern---is possible under the
condition $d^*\log(d/d^*)=o(n)$ as
$n\to\infty$~\cite{Wainwright09}. Furthermore, it is well known that
if $(d^*\log(d/d^*))/n$ remains bounded from
below by some positive constant when $n\to\infty$, then it is
impossible to consistently recover the sparsity
pattern~\cite{Verzelen10}. Thus, a tight condition exists
that describes in an exhaustive manner the interplay between the
quantities~$d^*$, $d$ and $n$ that guarantees the
existence of consistent estimators. The situation is very different in
the case of nonlinear regression, since, to
our knowledge, there is no result providing tight conditions for
consistent estimation of the sparsity pattern.

Lafferty and Wasserman~\cite{LaffertyWasserman} and Bertin and Lecu{\'
e}~\cite{BertinLecue}, in papers closely related to the present work,
considered the problem of
variable selection in nonparametric Gaussian regression model. They
proved the consistency of the proposed procedures under
some assumptions that---in the light of the present work---turn out to
be suboptimal. More precisely, Lafferty and Wasserman~\cite{LaffertyWasserman}
assumed the unknown regression function to be four times continuously
differentiable with bounded derivatives. The algorithm
they proposed, termed Rodeo, is a greedy procedure performing
simultaneously local bandwidth choice and variable selection.
Rodeo is shown to converge when the ambient dimension $d$ is $O({\log
n}/{\log\log n})$ while the intrinsic dimension $d^*$
does not increase with $n$. On the other hand, Bertin and Lecu{\'e}~\cite{BertinLecue}
proposed a procedure based on the $\ell_1$-penalization of
local polynomial estimators and proved its consistency when $d^*=O(1)$,
but $d$ is allowed to be as large as $\log n$, up to a
constant. They also had a weaker assumption on the regression function
merely assumed to belong to the Holder class
with smoothness $\beta>1$. To complete the picture, let us mention that
estimation and hypotheses testing problems for high-dimensional
nonparametric regression under sparse additive modeling were recently
addressed in~\cite{Koltchinskii10,Raskuttietal,GayraudIngster}.

This brief review of the literature reveals that there is an important
gap in consistency conditions for the linear regression
and for the nonlinear one. For instance, if the intrinsic dimension
$d^*$ is fixed, then the condition guaranteeing consistent
estimation of the sparsity pattern is $(\log d)/n\to0$ in linear
regression, whereas it is $d=O(\log n)$ in the nonparametric case.
While it is undeniable that the nonparametric regression is much more
complex than the linear one, it is, however, not easy to find
a justification to such an important gap between two conditions. The
situation is even worse in the case where $d^*\to\infty$.
In fact, for the linear model with at most polynomially increasing
ambient dimension $d=O(n^k)$, it is possible to estimate
the sparsity pattern for intrinsic dimensions $d^*$ as large as
$n^{1-\varepsilon}$, for some $\varepsilon>0$. In other words, the
sparsity index can be almost on the same order as the sample size. In
contrast, in nonparametric regression, there is no
procedure that is proved to converge to the true sparsity pattern when
both $n$ and $d^*$ tend to infinity, even if $d^*$
grows extremely slowly.

In the present work, we fill this gap by introducing a simple variable
selection procedure that selects the relevant variables
by comparing some quadratic functionals of empirical Fourier
coefficients to prescribed significance levels. Consistency of this
procedure is established under some conditions on the triplet
$(d^*,d,n)$, and the tightness of these conditions is proved. The
main take-away messages deduced from our results are the following:
\begin{itemize}
\item When the number of relevant variables $d^*$ is fixed
and the sample size $n$ tends to infinity, there exist
positive real numbers ${c}_*$ and ${c}^*$ such that (a) if $(\log
d)/n\le{c}_*$ the estimator proposed in Section~\ref{sec3}
is consistent and (b) no estimator of the sparsity pattern may be
consistent if $(\log d)/n\ge{c}^*$.
\item When the number of relevant variables $d^*$ tends to
infinity with $n\to\infty$, then there exist
real numbers $\underline{c}_i$ and $\bar c_i$, $i=1,2$ such that
$\underline{c}_1>0$, $\bar c_1>0$
and (a)~if $\underline{c}_1d^*+\log\log(d/d^*)-\log n<\underline{c}_2$
the estimator proposed in Section~\ref{sec3} is consistent
and (b)~no estimator of the sparsity pattern may be consistent if $\bar
c_1d^*+\log\log(d/d^*)- \log n >\bar c_2$.
\item In particular, if $d$ grows not faster than a
polynomial in $n$, then there exist
positive real numbers ${c}_0$ and ${c}^0$ such that (a) if $d^*\le
{c}_0\log n$, the estimator proposed in Section~\ref{sec3}
is consistent, and (b) no estimator of the sparsity pattern may be
consistent if $d^*\ge{c}^0\log n$.
\end{itemize}
In the regime of a growing intrinsic dimension $d^*\to\infty$ and a
moderately large ambient dimension $d=O(n^C)$, for some $C>0$,
we make a concentrated effort to get the constant $c_0$ as close as
possible to the constant $c^0$. This goal is reached for the model
of Gaussian white noise and, very surprisingly, it required from us to
apply some tools from complex analysis, such as the
Jacobi $\theta$-function and the saddle point method, in order to
evaluate the number of lattice points lying in a ball of an
Euclidean space with increasing dimension.

The rest of the paper is organized as follows. The notation and
assumptions necessary for stating our main results are presented in
Section~\ref{sec2}. In Section~\ref{sec3}, an estimator of the set of
relevant variables is introduced and its consistency
is established, in the case where the data come from the Gaussian white
noise model. The main condition required in the consistency
result involves the number of lattice points in a ball of a
high-dimensional Euclidean space. An asymptotic equivalent for this number
is presented in Section~\ref{sec4}. Results on impossibility of
consistent estimation of the sparsity pattern are derived in
Section~\ref{sec5}. Section~\ref{secadapt} is devoted to exploring
adaptation to the unknown parameters (smoothness and degree of
significance) and recovering minimax rates of separation. Then, in
Section~\ref{sec6}, we show that some of our results can be
extended to the model of nonparametric
regression. The relations between consistency and inconsistency results
are discussed in Section~\ref{sec7}. The technical parts of
the proofs are postponed to the \hyperref[app]{Appendix}.

\section{The problem formulation and the assumptions}\label{sec2}

We are interested in the variable selection task (also known as model
selection, feature selection, sparsity pattern estimation)
in the context of high-dimensional nonlinear regression. Let $\sff
\dvtx [0,1]^d\to\RR$ denote the unknown regression function. We assume
that the number of variables $d$ is very large, possibly much larger
than the sample size $n$, but only a small
number of these variables contribute to the fluctuations of the
regression function $\sff$.

To be more precise, we assume that for some small subset $J$ of the
index set $\{1,\ldots,d\}$
satisfying $\operatorname{Card}(J)\le d^*$, there is a function $\bar
\sff\dvtx \RR^{\mathrm{Card}(J)}\to\RR$ such that
\[
\sff(\mathbf{x})=\bar\sff(\mathbf{x}_J) \qquad\forall\mathbf{x}
\in\RR^d,
\]
where $\mathbf{x}_J$ stands for the subvector of $\mathbf{x}$
obtained by removing
from $\mathbf{x}$ all the coordinates with indices lying
outside $J$. In what follows, we allow $d$ and $d^*$ to depend on $n$,
but we will not always indicate this dependence
in notation. Note also that the genuine intrinsic dimension is
$\operatorname{Card}(J)$; $d^*$ is merely a known upper bound on
the intrinsic dimension. In what follows, we use the standard notation
for the vector and sequence norms:
\begin{eqnarray*}
\|\mathbf{x}\|_0&=&\sum_{j}
\mathbf1(x_j\not=0),\qquad \|\mathbf{x}\|_p^p=
\sum_{j} |x_j|^p\qquad \forall p
\in[1,\infty), \\
\|\mathbf{x}\|_\infty&=&\sup_{j}
|x_j|
\end{eqnarray*}
for every $\mathbf{x}\in\RR^d$ or $\mathbf{x}\in\RR^\NN$.

Let us stress right away that the primary aim of this work is to
understand when it is possible to estimate the sparsity pattern $J$
(with theoretical guarantees on the convergence of the estimator) and
when it is impossible. The estimator that we will define in the
next sections is intended to show the possibility of consistent
estimation, rather than to provide a practical procedure for recovering
the sparsity pattern. Therefore, the estimator will be allowed to
depend on different constants appearing in conditions imposed on the
regression function $\sff$ and on some characteristics of the noise.

To make the consistent estimation of the set $J$ realizable, we impose
some smoothness and identifiability assumptions on $\sff$.
In order to describe the smoothness assumption imposed on $\sff$, let
us introduce the trigonometric Fourier basis, $\varphi_{\mathbf{0}}\equiv1$ and
%
\begin{equation}
\label{trig} \varphi_\bk(\mathbf{x})= \cases{ \sqrt{2}\cos(2\pi\bk
\cdot\mathbf{x}), &\quad $\bk\in\bigl(\ZZ^d\bigr)_+,$\vspace*{2pt}
\cr
\sqrt{2}\sin(2\pi\bk\cdot\mathbf{x}), &\quad $-\bk\in\bigl(\ZZ^d\bigr)_+
,$}
\end{equation}
where $(\ZZ^d)_+$ denotes the set of all $\bk\in\ZZ^d\setminus\{0\}$
such that the first nonzero element of $\bk$ is positive, and
$\bk\cdot\mathbf{x}$ stands for the usual inner product in $\RR^d$. In what
follows, we use the notation $\langle\cdot,\cdot\rangle$
for designing the scalar product in $L^2([0,1]^d;\RR)$, that is,
$\langle\sfh,\tilde\sfh\rangle=\int_{[0,1]^d} \sfh(\mathbf{x})\tilde\sfh(\mathbf{x}
) \,d\mathbf{x}$
for every $\sfh,\tilde\sfh\in L^2([0,1]^d;\RR)$. Using this
orthonormal Fourier basis, we define
\[
\Sigma_L= \biggl\{f\dvtx \sum_{\bk\in\ZZ^d}
k_j^{2}\langle\sff ,\varphi_\bk
\rangle^2\le L; \forall j\in\{1,\ldots,d\} \biggr\}.
\]
To ease notation, we set $\theta_\bk[\sff]=\langle\sff,\varphi_\bk
\rangle$ for all $\bk\in\ZZ^d$. In addition to the smoothness, we
need also to require that the relevant variables are sufficiently
relevant for making their identification possible. This is done
by means of the following condition.
\begin{longlist}[{[C1$(\kappa,L)$]}]
\item[{[C1$(\kappa,L)$]}] The regression function $\sff$ belongs to
$\Sigma_L$. Furthermore, for some subset $J\subset\{1,\ldots,d\}$ of
cardinality $\le d^*$, there exists a function\vadjust{\goodbreak} $\bar\sff\dvtx \RR^{\mathrm{Card}(J)}\to\RR$ such that
$\sff(\mathbf{x}) =\bar\sff(\mathbf{x}_J)$, $\forall
\mathbf{x}\in\RR^d$, and it holds that
%
\begin{equation}
\label{ident} Q_j[\sff]=\sum_{\bk: k_j\neq0}
\theta_\bk[\sff]^2\geq\kappa \qquad\forall j\in J.
\end{equation}
\end{longlist}
One easily checks that $Q_j[\sff]=0$ for every $j$ that does not lie in
the sparsity pattern. This provides a characterization of the
sparsity pattern as the set of indices of nonzero coefficients of the
vector $\bQ[\sff]=(Q_1[\sff],\ldots,Q_d[\sff])$.

Prior to describing the procedures for estimating $J$, let us comment
on condition~{[C1]}. It is important to note that
the identifiability assumption (\ref{ident}) can be rewritten as $\int_{[0,1]^d} (\sff(\mathbf{x})-\int_0^1\sff(\mathbf{x})
\,dx_j )^2
\,d\mathbf{x}\ge\kappa$ and, therefore, is not intrinsically related
to the
basis we have chosen. In the case of continuously differentiable
and $1$-periodic function $\sff$, the smoothness assumption $\sff\in
\Sigma_L$ as well can be rewritten without using the trigonometric
basis, since $\sum_{\bk\in\ZZ^d} k_j^{2}\theta_\bk[\sff]^2=(2\pi
)^{-2}\int_{[0,1]^d}[\partial_j \sff(\mathbf{x})]^2 \,d\mathbf{x}$. Thus
condition {[C1]} is essentially a constraint on the function
$\sff$ itself and not on its representation in the specific
basis of trigonometric functions.

The results of this work can be extended with minor modifications to
other types of smoothness conditions imposed on $\sff$,
such as H\"older continuity or Besov-regularity. In these cases the
trigonometric basis (\ref{trig}) should be replaced by a basis
adapted to the smoothness condition (spline, wavelet, etc.).
Furthermore, even in the case of Sobolev smoothness, one can replace
the set $\Sigma_L$ corresponding to smoothness order $1$ by any Sobolev
ellipsoid of smoothness $\beta>0$; see, for instance, \cite
{LaetCRASS} where
the case $\beta=2$ is explored. Roughly speaking, the role of the
smoothness assumption is to reduce the statistical model with
infinite-dimensional
parameter $\sff$ to a finite-dimensional model having good
approximation properties. Any value of smoothness order $\beta>0$ leads
to this reduction. The value $\beta=1$ is chosen for simplicity of
exposition only.

\section{Idealized setup: Gaussian white noise model}\label{sec3}

To convey the main ideas without taking care of some technical details,
we start by focusing our attention on
the Gaussian white noise model that was proved to be asymptotically
equivalent to the model of regression
\cite{BrownLow96,Reiss08}, as well as to other nonparametric models~\cite{Brown04,DalReiss06}. Thus, we assume that
the available data consists of the Gaussian process $\{\bY(\phi)\dvtx \phi
\in L^2([0,1]^d;\RR)\}$
such that
\[
\E_\sff\bigl[\bY(\phi)\bigr]=\int_{[0,1]^d} \sff(
\mathbf{x}) \phi (\mathbf{x}) \,d\mathbf{x},\qquad \operatorname{Cov}_f
\bigl(\bY(\phi),\bY\bigl(\phi'\bigr)\bigr)=\frac1n\int
_{[0,1]^d} \phi (\mathbf{x})\phi'(\mathbf{x}) \,d
\mathbf{x}.
\]
It is well known that these two properties uniquely characterize the
probability distribution of a Gaussian process. An alternative
representation of $\bY$ is
\[
dY(\mathbf{x})=\sff(\mathbf{x}) \,d\mathbf{x}+n^{-1/2}\,dW(
\mathbf{x}),\qquad \mathbf{x}\in[0,1]^d,
\]
where $W(\mathbf{x})$ is a $d$-parameter Brownian sheet. Note that minimax
estimation and detection of the function $\sff$ in this set-up
(but without sparsity assumption) was studied by~\cite{Ingster07}.

\subsection{Estimation of $J$ by multiple hypotheses testing}
We intend to tackle the variable selection problem by multiple
hypotheses testing; each hypothesis
concerns a group of the Fourier coefficients of the observed signal and
suggests that all the elements
within the group are zero. The rationale behind this approach is the
following simple observation:
since the trigonometric basis is orthonormal and contains the constant function,
%
\begin{equation}
\label{crit1} j\notin J \quad\Longleftrightarrow\quad \theta_\bk[\sff]=\langle
f,\varphi_\bk \rangle=0 \qquad\forall\bk \mbox{ s.t. } k_j
\not=0.
\end{equation}
This observation entails that if the intrinsic dimension $|J|$ is small
as compared to~$d$, then the sequence
of Fourier coefficients is sparse. Furthermore, as explained below,
there is a sort of group sparsity with overlapping groups.

For every $\ell\in\{1,\ldots,d^*\}$, we denote by $P_\ell^d$ the
set of
all subsets $I$ of $\{1,\ldots,d\}$
having exactly $\ell$ elements: $P_\ell^d= \{I\subset\{1,\ldots
,d\}
\dvtx \operatorname{Card}(I)=\ell \}$.
For every multi-index $\bk\in\ZZ^d$, we denote by $\supp(\bk)$ the set
of indices corresponding to nonzero
entries of $\bk$. To define the blocks of coefficients $\theta_{\bk}$
that will be tested for significance,
we introduce the following notation: for every $I\subset\{1,\ldots,d\}
$ and for every $j\in I$, we set
\[
V_{I}^j[\sff]= \bigl(\theta_\bk[\sff]\dvtx j
\in\supp(\bk)\subset I \bigr).
\]
It follows from (\ref{crit1}) that the characterization
%
\begin{equation}
\label{crit2} j\notin J \quad\Longleftrightarrow\quad \max_I
\bigl\|V_I^j[\sff] \bigr\|_p=0,
\end{equation}
holds true for every $p\in[0,+\infty]$. Furthermore, again in view of
(\ref{crit1}), the maximum over $I$ of the
norms $ \|V_I^j[\sff] \|_p$ is attained when $I=J$ and is equal
to the maximum over
all subsets $I$ such that $\operatorname{Card}(I)\le d^*$. Summarizing
these arguments, we can formulate the problem
of variable selection as a problem of testing $d$ null hypotheses
%
\begin{equation}
\label{H0} H_{0j} \dvtx \bigl\|V_I^j[\sff]
\bigr\|_p=0 \qquad \forall I\subset\{1,\ldots,d\} \mbox{ such that }
\operatorname{Card}(I)\le d^*.
\end{equation}
If the hypothesis $H_{0j}$ is rejected, then the $j$th covariate is
declared as relevant. Note that by virtue of
assumption [{C1}], the alternatives can be written as
%
\begin{equation}
\label{H1} H_{1j} \dvtx \bigl\|V_I^j[\sff]
\bigr\|_2^2\ge\kappa\qquad\mbox{for some } I\subset\{1,\ldots,d\}
\mbox{ such that } \operatorname{Card}(I)\le d^*.\hspace*{-35pt}
\end{equation}
Our estimator is based on this characterization of the sparsity
pattern. If we denote by $y_\bk$ the
observable random variable $\bY(\varphi_\bk)$, we have
%
\begin{equation}
\label{WGN} y_\bk=\theta_\bk[\sff]+n^{-1/2}
\xi_\bk,\qquad \theta_\bk=\langle f,\varphi_\bk
\rangle, \bk\in\ZZ^d,
\end{equation}
where $\{\xi_\bk;\bk\in\ZZ^d\}$ form a countable family of independent
Gaussian random variables with zero mean and
variance equal to one. According to this property, $y_\bk$ is a good
estimate of $\theta_\bk[\sff]$: it is unbiased
and with a mean squared error equal to $1/n$. Using the plug-in
argument, this suggests to estimate $V_I^j$ by
$\widehat V_I^j=(y_\bk\dvtx j\in\supp(\bk)\subset I)$ and the norm of
$V_I^j$ by the norm of $\widehat V_I^j$. However,
since this amounts to estimating an infinite-dimensional vector, the
error of estimation will be infinitely large.
To cope with this issue, we restrict the set of indices for which
$\theta_\bk$ is estimated by $y_\bk$ to a finite set,
outside of which $\theta_\bk$ will be merely estimated by $0$. Such a
restriction is justified by the fact that
$\sff$ is assumed to be smooth: Fourier coefficients corresponding to
very high frequencies are very small.

Let us fix an integer $m>0$, the cut-off level, and denote, for $j\in
I\subset\{1,\ldots, d\}$,
\[
S^j_{m,I}= \bigl\{\bk\in\ZZ^d \dvtx \|\bk
\|_{2}\le m \mbox{ and } \{ j\} \subset\operatorname{supp} (\bk)\subset
I \bigr\}.
\]
Since the alternatives $H_{1j}$ are concerned with the 2-norm, we build
our test statistic on an estimate of
the norm $\|V_I^j[\sff]\|_2$. To this end, we introduce
\[
\widehat{Q}_{m,I}^j=\sum_{\bk\in S^j_{m,I}}
\biggl(y_\bk^2-\frac
{1}{n} \biggr),
\]
which is an unbiased estimator\vspace*{-2pt} of $Q_{m,I}^{j}=\sum_{\bk\in
S^j_{m,I}}\theta_\bk^2$. Note that when
$m\to\infty$, the quantity $Q_{m,I}^{j}$ approaches\vspace*{-2pt} $\|V_I^j[\sff]\|_2^2$. It is clear that larger values of $m$
lead to a smaller bias while the variance get increased. Moreover, the
variance of $\widehat{Q}_{m,I}^j$
is proportional to the cardinality of the set $S^j_{m,I}$. The latter
is an increasing function of $\operatorname{Card}(I)$.
Therefore, if we aim at getting comparable estimation accuracies when
estimating the functionals
$\|V_I^j[\sff]\|_2^2$ by $\widehat{Q}_{m,I}^j$ for various $I$'s, it is
reasonable to make the cut-off level $m$ vary
with the cardinality of $I$.

Thus, we consider a multivariate cut-off $\mm=(m_1,\ldots,m_{d^*})\in
\NN^{d^*}$. For a subset~$I$ of cardinality $\ell\le d^*$,
we test significance of the vector $V_I^j[\sff]$ by comparing its
estimate $\widehat{Q}_{m_\ell,I}^j$ with a prescribed
threshold $\lambda_\ell$. This leads us to define an estimator of the
set $J$ by
\[
\widehat{J}_{n}(\mm,\llambda)= \Bigl\{j\in\{1,\ldots, d\} \dvtx
\max_{\ell
\leq d^*} \lambda_\ell^{-1} \max_{I\in P^d_{\ell}}
\widehat{Q}_{m_\ell,I}^j\ge1 \Bigr\},
\]
where $\mm=(m_1,\ldots,m_{d^*})\in\NN^{d^*}$ and $\llambda
=(\lambda_1,\ldots,\lambda_{d^*})\in\RR^{d^*}_+$ are
two vectors of tuning parameters. As already mentioned, the role of
$\mm
$ is to ensure that the truncated sums
$Q_{m,J}^j$ do not deviate too much from the complete sums $Q_J^j$.
Quantitatively speaking, for a given $\tau>0$,
we would like to choose $m_\ell$'s so that $Q_{m_s,J}^j\ge\kappa\tau
/\tau+1$, where $s=\operatorname{Card}(J)$. This
guarantee can be achieved due to the smoothness assumption. Indeed, as
proved in (\ref{minor}) (cf. Appendix~\ref{appA1}),
it holds that
\begin{eqnarray*}
Q_{m_s,J}^{j}\geq\kappa- m_s^{-2} Ls\qquad
\forall j\in J.
\end{eqnarray*}
Therefore, choosing $m_\ell= (\ell L(1+\tau) /\kappa )^{1/2}$,
for every $\ell=1,\ldots,d^*$, entails the inequality
$Q_{m_s,J}^j\ge\kappa\tau/\tau+1$, which indicates that the relevance
of variables is not affected too much
by the truncation.

Pushing further the analogy with the hypotheses testing, we define type
I error of an estimator $\widehat J_n$ of $J$
as the one of having $\widehat{J}_n\not\subset J$, that is, classifying
some irrelevant variables as relevant.
The type II error is then that of having $J\not\subset\widehat J$,
which amounts to classifying some relevant variables
as irrelevant. As in the testing problem, handling the type I error is
easier since the distribution of the test statistic
is independent of $\sff$. In fact, this is the max of a finite family
of random variables drawn from translated and
scaled $\chi^2$-distributions. Using the Bonferroni adjustment leads to
the following control of the first kind error.

\begin{proposition}\label{prop0}
Let us denote by $N(\ell, \gamma)$ the cardinality of the set $\{\bk
\in
\ZZ^{\ell}\dvtx \|\bk\|_2^2\le\gamma\ell\ \&\ k_1\not=0\}$.
If for some $A>1$ and for every $\ell=1,\ldots,d^*$,
%
\begin{equation}
\label{lambda} \lambda_\ell\ge\frac{2\sqrt{A N(\ell,m_\ell^2/\ell)d^*\log
(2ed/d^*)}+2Ad^*\log(2ed/d^*)}{n},
\end{equation}
then the type I error $\Pb (\widehat{J}_n(\mm,\llambda)\not
\subset
J )$ is upper-bounded by $(2ed/d^*)^{-d^*(A-1)}$, and
therefore tends to 0 as $d\to+\infty$.
\end{proposition}

This proposition shows that the type I error of a variable selection
procedure may be made small by choosing a sufficiently
high threshold. By doing this, we run the risk to reject $H_{0j}$ very
often and to drastically underestimate the set of relevant
variables. The next result establishes a necessary condition, which
will be shown to be tight, ensuring that such an
underestimation does not occur.

\begin{theorem}\label{prop2}
Let condition \textup{[C1$(\kappa,L)$]} be satisfied with some known
constants $\kappa>0$ and $L<\infty$, and let $s=\operatorname{Card}(J)$.
For some real numbers $\tau>0$ and $A>1$, set
$m_\ell= (\ell L(1+\tau) /\kappa )^{1/2}$, $\ell=1,\ldots,d^*$,
and define $\lambda_\ell$ to be equal to the right-hand side of (\ref
{lambda}).
If the condition
%
\begin{equation}
\label{ordre2} 4\lambda_s\le{\kappa\tau}/{(1+\tau)}
\end{equation}
is fulfilled, then $\widehat J_{n}(\mm,\llambda)$ is consistent and
satisfies the inequalities $\Pb (\widehat J_{n}(\mm,\break\llambda
)\not
\supset J )\leq2(2ed/d^*)^{-d^*(A-1)}$ and
$\Pb (\widehat J_{n}(\mm,\llambda)\not=J )\leq3
(2ed/d^*)^{-d^*(A-1)}$.
\end{theorem}

Condition (\ref{ordre2}), ensuring the consistency of the variable
selection procedure $\widehat J_n$,
admits a very natural interpretation: It is possible to detect relevant
variables if the degree of
relevance $\kappa$ is larger than a multiple of the threshold $\lambda_s$, the latter being chosen
according to the noise level.

A first observation is that this theorem provides interesting insight
to the possibility of consistent
recovery of the sparsity pattern $J$ in the context of fixed intrinsic
dimension. In fact, when $d^*$
remains bounded from above when $n\!\to\!\infty$ and $d\!\to\!\infty$, then we
get that
$\Pb(\widehat J_1(m,\lambda)= J)\!\to\!_{n,d\to\infty} 1$ provided
that\looseness=-1
%
\begin{equation}
\label{eq9a} \log d\le\operatorname{Const}\cdot n.
\end{equation}\looseness=0
Although we did not find (exactly) this result in the statistical
literature on variable selection, it can be
checked that (\ref{eq9a}) is a necessary and sufficient condition for
recovering the sparsity pattern $J$
in linear regression with fixed sparsity $d^*$ and growing dimension
$d$ and sample size $n$. Thus, in the
regime of fixed or bounded~$d^*$, the sparsity pattern estimation in
nonparametric regression is not
more difficult than in the parametric linear regression, as far as only
the consistency of estimation
is considered and the precise value of the constant in (\ref{eq9a}) is
neglected. Furthermore, there is
a simple estimator $\widehat J_{n}^{(1)}$ of $J$ (cf. equation (3) in
\cite{LaetCRASS}), which is provably
consistent under condition (\ref{eq9a}). This estimator can be seen as
a procedure of testing hypotheses
$H_{0j}$ of form (\ref{H0}) with $p=\infty$ and, therefore, it does not
really exploit the structure of
the Fourier coefficients of the regression function. To some extent,
this is the reason why in the regime
of growing intrinsic dimension $d^*\to\infty$, the estimator
$\widehat
J_{n}^{(1)}$ proposed by
\cite{LaetCRASS} is no longer optimal.

In fact, when $d^*\to\infty$, the term $N(s,m^2_s/s)$ present in
(\ref
{ordre2}) tends to infinity as well.
Furthermore, as we show in Section~\ref{sec4}, this convergence takes
place at an exponential rate in $d^*$.
Therefore, in this asymptotic set-up it is crucial to have the right
order of $N(s,m^2_s/s)$ in the condition
that ensures the consistency. As shown in Section~\ref{sec5}, this is
the case for condition (\ref{ordre2}).

\begin{remark}
An apparent drawback of the estimator $\widehat J_n$ is the large
dimensionality of tuning parameters involved
in $\widehat J_{n}$. However, Theorem~\ref{prop2} reveals that for
achieving good selection power, it is
sufficient to select the $2d^*$-dimensional tuning parameter $(\mm
,\llambda)$ on a one-dimensional curve
parameterized by $\vartheta=L(1+\tau)/\kappa$. Indeed, once the value
of $\vartheta$ is given,
Theorem~\ref{prop2} advocates for choosing
%
\begin{eqnarray}
\label{mlambda} m_\ell&=&(\ell\vartheta)^{1/2} \quad\mbox{and}
\nonumber
\\[-8pt]
\\[-8pt]
\nonumber
\lambda_\ell&=&\frac
{2\sqrt{A
N(\ell,\vartheta)d^*\log(2ed/d^*)}+2Ad^*\log(2ed/d^*)}{n}
\end{eqnarray}
for every $\ell=1,\ldots,d^*$. As discussed in Section \ref
{ssecadapt}, this property allows us to relax the requirement
that the values $L$ and $\kappa$ involved in [C1] are known
in advance.\vadjust{\goodbreak}
\end{remark}

\begin{remark}
The result of the last theorem is in some sense adaptive w.r.t. the
unknown sparsity. Indeed, while the estimator $\widehat J_{n}$
involves $d^*$, which is merely a known upper bound on the true
sparsity $s=\operatorname{Card}(J)$ and may be significantly larger than $s$,
it is the true sparsity $s$ that appears in condition (\ref{ordre2}) as
a first argument of the quantity $N(\cdot, \vartheta)$.
This point is important given the exponential rate of divergence of
$N(\cdot,\vartheta)$ when its first argument tends
to infinity. On the other hand, if condition (\ref{ordre2}) is
satisfied with $N(d^*,\vartheta)$ instead of $N(\operatorname
{Card}(J),\vartheta)$,
then the consistent estimation of $J$ can be achieved by a slightly
simpler procedure,
\[
\widetilde{J}_{n}(\mm,\llambda)= \Bigl\{j\in\{1,\ldots, d\} \dvtx
\max_{I\in P^d_{d^*}} \widehat{Q}_{m_{d^*},I}^j\ge
\lambda_{d^*} \Bigr\}.
\]
The proof of this statement is similar to that of Theorem~\ref{prop2}
and will be omitted.
\end{remark}

\section{Counting lattice points in a ball}\label{sec4}

\begin{figure}

\includegraphics{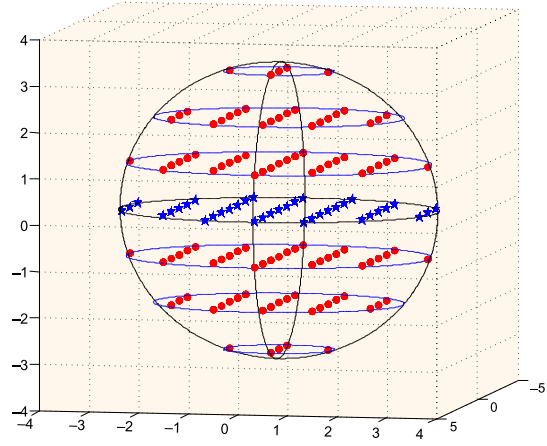}

\caption{Lattice points in a ball of radius $R=\gamma d^*=3.2$ in the
three dimensional space ($d^*=3$). Red points are those
of $\mathcal C_1(d^*,\gamma)\setminus\mathcal C_2(d^*,\gamma)$ while
blue stars are those of $\mathcal C_2(d^*,\gamma)$. In this
example, $N(d^*,\gamma)=N(3,1.07)=110$.}\label{fig01}
\end{figure}

The aim of the present section is to investigate the properties of the
quantity $N(d^*,\gamma)$ that is
involved in the conditions ensuring the consistency of the proposed
procedures. Quite surprisingly, the
asymptotic behavior of $N(d^*,\gamma)$ turns out to be related to the
Jacobi $\theta$-function. To
show this, let us introduce some notation. For a positive number
$\gamma
$, we set
\begin{eqnarray*}
\mathcal{C}_1\bigl(d^*,\gamma\bigr)&=& \bigl\{\bk\in\ZZ^{d^*}
\dvtx k_1^{2}+\cdots+k_{d^*}^{2}\leq
\gamma d^* \bigr\}, \\ \mathcal{C}_2\bigl(d^*,\gamma\bigr)&=& \bigl\{\bk\in
\mathcal{C}_1\bigl(d^*,\gamma\bigr) \dvtx k_1=0 \bigr\}
\end{eqnarray*}
along with $N_1(d^*,\gamma)=\operatorname{Card}\mathcal
{C}_1(d^*,\gamma
)$ and $N_2(d^*,\gamma)=\operatorname{Card}\mathcal{C}_2(d^*,\gamma)$.
In simple words, $N_1(d^*,\gamma)$ is the number of lattice points
lying in the $d^*$-dimensional ball
with radius $(\gamma d^*)^{1/2}$ and centered at the origin, while
$N_2(d^*,\gamma)$ is the number of (integer) lattice
points lying in the $(d^*-1)$-dimensional ball with radius $(\gamma
d^*)^{1/2}$ and centered at the origin; see Figure~\ref{fig01} for an illustration. With this
notation, the quantity $N(\ell,\cdot)$ of Theorem~\ref{prop2} can be
written as
$N_1(\ell,\cdot)-N_2(\ell,\cdot)$. By volumetric arguments, one can
check that
$V(d^*) (\sqrt{\gamma}-1)^{d^*} (d^*)^{d^*/2}\le N_1(d^*,\gamma)\le
V(d^*) (\sqrt{\gamma}+1)^{d^*} (d^*)^{d^*/2}$, where
$V(d^*)=\pi^{d^*/2}/\Gamma(1+d^*/2)$ is the volume of the unit ball in
$\RR^{d^*}$. Furthermore, similar bounds hold true
for $N_2(d^*,\gamma)$ as well. Unfortunately, when $d^*\to\infty$,
these inequalities are not accurate enough to yield
nontrivial results in the problem of variable selection we are dealing
with. This is especially true for the results
on impossibility of consistent estimation stated in Section~\ref{sec5}.

In order to determine the asymptotic behavior of $N_1(d^*,\gamma)$ and
$N_2(d^*,\gamma)$ when $d^*$ tends to infinity,
we will rely on their integral representation through Jacobi's $\theta
$-function. Recall that the latter is given
by $\sfh(z)=\sum_{r\in\mathbb{Z}}z^{r^{2}}$, which is well defined for
any complex number $z$ belonging to the unit ball
$|z|<1$. To briefly explain where the relation between $N_i(d^*,\gamma
)$ and the $\theta$-function comes from, let
us denote by $\{a_r\}$ the sequence of coefficients of the power series
of $\sfh(z)^{d^*}$, that is, $\sfh(z)^{d^*}=
\sum_{r\geq0} a_r z^r$. One easily checks that $\forall r\in\mathbb{N}$,
$a_r=\operatorname{Card}\{\bk\in\mathbb{Z}^{d^*} \dvtx k_1^{2}+\cdots+k_{d^*}^{2}=r\}$. Thus, for every $\gamma$ such
that $\gamma d^*$ is integer, we have $N_1(d^*,\gamma)=\sum_{r=0}^{\gamma d^*} a_r$. As a consequence of Cauchy's theorem,
we get
\[
N_1\bigl(d^*,\gamma\bigr)=\frac{1}{2\pi i}\oint\frac{\sfh
(z)^{d^*}}{z^{\gamma
d^*}}
\frac{dz}{z(1-z)},
\]
where the integral is taken over any circle $|z|=w$ with $0<w<1$.
Exploiting this representation and applying the saddle-point
method thoroughly described in~\cite{Dieudonne}, we get the following result.

\begin{figure}

\includegraphics{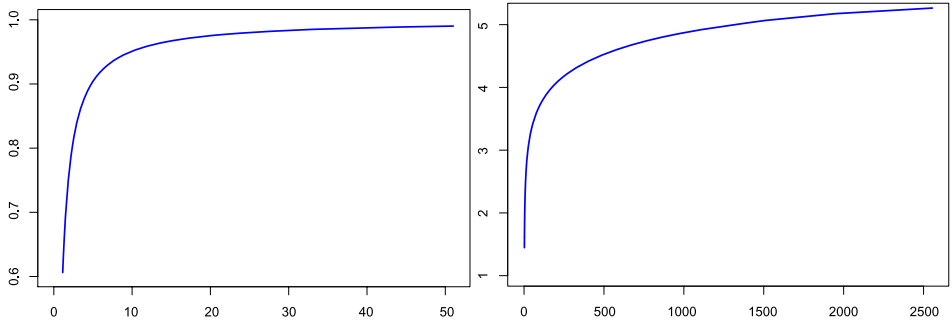}

\caption{The plots of mappings $\gamma\mapsto z_\gamma$ and $\gamma
\mapsto\sfl_\gamma(\zb)$. One can observe that both functions are
increasing,
the first one converges to $1$ very rapidly, while the second one seems
to diverge very slowly. }
\label{fig1}
\end{figure}

\begin{proposition}\label{prop1}
Let $\gamma>0$ be an integer and let $\sfl_\gamma(z)=\log{\sfh
(z)}-\gamma\log z$.

\begin{longlist}[(1)]
\item[(1)] There is a unique solution $\zb$ in $(0,1)$ to the equation
$\sfl_\gamma'(z)=0$. Furthermore,
the function $\gamma\mapsto z_\gamma$ is increasing and $\sfl_\gamma''(z)>0$.

\item[(2)] For $i=1,2$, the following equivalences hold true:
\[
N_i\bigl(d^*,\gamma\bigr)= \biggl(\frac{\sfh(z_\gamma)}{z_\gamma^\gamma
}
\biggr)^{d^*}\frac{1+o(1)}{\sfh(z_\gamma)^{i-1}\zb(1-\zb)(2\sfl_\gamma''(\zb
)\pi d^*)^{1/2}}
\]
as $d^*$ tends to infinity.
\end{longlist}
\end{proposition}
Hereafter, it will be useful to note that the second part of
Proposition~\ref{prop1} yields
%
\begin{eqnarray}
\label{eqlog} \log \bigl(N_1\bigl(d^*,\gamma\bigr)-N_2
\bigl(d^*,\gamma\bigr) \bigr)  = d^*\sfl_\gamma (z_\gamma)-
\tfrac12\log d^*+c_\gamma+o(1)
\nonumber
\\[-8pt]
\\[-8pt]
\eqntext{\mbox{as }d^*\to\infty,}
\end{eqnarray}
with $c_\gamma=\log (\frac{\sfh(\zb)-1}{\sfh(\zb)\zb(1-\zb
)\sqrt {2\pi\sfl''_\gamma(\zb)}} )$.
Furthermore, while the asymptotic equivalences of Proposition \ref
{prop1} are established for integer values of $\gamma>0$,
relation $\log (N_1(d^*,\gamma)-N_2(d^*,\gamma) ) = d^*\sfl_\gamma(z_\gamma)(1+o(1))$ holds true for any positive real number
$\gamma$~\cite{Mazo}. In order to get an idea of how the terms $\zb$ and
$\sfl_\gamma(\zb)$ depend on~$\gamma$, we depicted in Figure \ref
{fig1} the plots of these quantities as functions of $\gamma>0$.

Combining relation (\ref{eqlog}) with Theorem~\ref{prop2}, we get the
following result.
%
\begin{corollary}\label{cor2}
Let condition \textup{[C1$(\kappa,L)$]} be satisfied with some known
constants $\kappa>0$ and $L<\infty$. Consider the asymptotic
set-up in which both $d=d_n$ and $d^*=d^*_n$ tend to infinity as $n\to
\infty$. Assume that $d$ grows at a sub-exponential rate
in $n$, that is, $\log\log d=o(\log n)$. If
\[
\limsup_{n\to\infty} \frac{d^*}{\log n} <\frac{2}{\sfl_\gamma
(z_\gamma)}
\]
with $\gamma=L/\kappa$, then consistent estimation of $J$ is possible
and can be achieved, for instance, by the estimator
$\widehat J_{n}$.
\end{corollary}

\section{Tightness of the assumptions}\label{sec5}

In this section, we focus our attention on the functional class $\Sigma
(\kappa,L)$ of all functions satisfying
assumption [C1($\kappa,L$)]. For emphasizing that $J$ is the
sparsity pattern of the function
$\sff$, we write $J_\sff$ instead of $J$. We assume that
$s=\operatorname{Card}(J)=d^*$. The goal is to
provide conditions under which the consistent estimation of the
sparsity support is impossible, that is,
there exists a constant $c>0$ and an integer $n_0\in\NN$ such that, if
$n\ge n_0$,
\[
\inf_{\widetilde{J}}\sup_{\sff\in\Sigma(\kappa,L)} \Pb_\sff (\widetilde {J}\neq
J_\sff)\geq c,
\]
where the $\inf$ is over all possible estimators of $J_\sff$. To this
end, we
introduce a set of $M+1$ probability distributions $\mu_0,\ldots,\mu_{M}$ on $\Sigma(\kappa,L)$ and use the fact that
%
\begin{equation}
\label{minor1} \inf_{\widetilde{J}}\sup_{\sff\in\widetilde\Sigma(\kappa,L)} \Pb_\sff (
\widetilde{J}\neq J_\sff) \ge\inf_{\widetilde{J}}\frac1{M}\sum
_{\ell=1}^M\int_{\Sigma
(\kappa,L)}
\Pb_\sff(\widetilde{J}\neq J_\sff) \mu_\ell(d
\sff).
\end{equation}
These measures $\mu_\ell$ will be chosen in such a way that for each
$\ell\ge1$ there is a set $J_\ell$ of cardinality $d^*$
such that $\mu_\ell\{J_\sff=J_\ell\}=1$ and all the sets
$J_1,\ldots
,J_M$ are distinct. The measure $\mu_0$ is the Dirac measure in
$0$. Considering these $\mu_\ell$s as ``priors'' on $\Sigma(\kappa,L)$
and defining the corresponding ``posteriors''
$\PP_0,\PP_1,\ldots,\PP_M$ by
\[
\PP_\ell(A)=\int_{\Sigma(\kappa,L)} \Pb_\sff(A)
\mu_\ell(d\sff )\qquad \mbox{for every measurable set } A\subset
\RR^n,
\]
we can write inequality (\ref{minor1}) as
%
\begin{equation}
\label{minor2} \inf_{\widetilde{J}}\sup_{\sff\in\Sigma(\kappa,L)} \Pb_\sff (
\widetilde {J}\neq J_\sff) \ge\inf_{\psi}\frac1{M}\sum
_{\ell=1}^M \PP_\ell(\psi\not=
\ell),
\end{equation}
where the $\inf$ is taken over all random $\psi$ taking values in $\{
0,\ldots,M\}$. The latter $\inf$ will be controlled using
a suitable version of the Fano lemma. To state it, we denote by $\KL(P,
Q)$ the Kullback--Leibler divergence between two
probability measures~$P$ and $Q$ defined on the same probability space.

\begin{lemma}[(Corollary 2.6 of~\cite{Tsybakov09})]\label{lemfano}
Let $M\ge3$ be an integer, $(\mathcal X,\mathcal A)$ be a measurable
space and let $P_0,\ldots,P_M$ be probability measures on
$(\mathcal X,\mathcal A)$. Let us set $\bar{p}_{e,M}=\inf_\psi
M^{-1}\sum_{\ell=1}^M P_\ell (\psi\not= \ell )$,
where the $\inf$ is taken over all measurable functions $\psi
\dvtx \mathcal
X\to \{1,\ldots,M \}$. If for some $0<\alpha< 1$,
$\frac{1}{M+1}\sum_{\ell=1}^{M}\KL (P_\ell,\break P_0 )\leq
\alpha\log
M$, then $\bar{p}_{e,M}\ge\frac12-\alpha$.
\end{lemma}

We apply this lemma with $\mathcal X$ being the set of all arrays
$ \mathbf{y}
=\{y_\bk\dvtx \bk\in\ZZ^d\}$ such that
for some $K>0$ the entries $y_\bk=0$ for every $\bk$ larger than $K$ in
$\ell_2$-norm. It follows from Fano's
lemma that one can deduce a lower bound on $\bar p_{e,M}$, the quantity
we are interested in, from an
upper bound on the average Kullback--Leibler divergence between $\PP_\ell$ and $\PP_0$. With these tools at
hand, we are in a position to state the main result on the
impossibility of consistent estimation
of the sparsity pattern in the case when the conditions of Theorem \ref
{prop2} are violated.

\begin{theorem}\label{thm3}
Assume that $\vartheta=L/\kappa>1$ and $({{d}\atop{d^*}})\ge3$. Let
$\gammaL$ be the largest integer satisfying $\gamma (1+ (\sfh
(z_\gamma)-1 )^{-1} )\le\vartheta$,
where the Jacobi $\theta$-function $\sfh$ and $z_\gamma$ are those
defined in Section~\ref{sec4}.
\begin{longlist}[(ii)]
\item[(i)] If for some $\alpha\in(0,1/2)$,
%
\begin{equation}
\label{ordre3} \frac{N(d^*,\gammaL)d^*\log(d/d^*)}{n^2}\geq\frac{\vartheta
}{\alpha
\gammaL} \kappa^2,
\end{equation}
then, for $d^*$ large enough,
$\inf_{\widetilde{J}}\sup_{\sff\in\Sigma} \Pb_\sff
(\widetilde
{J}\neq J_\sff )\geq\frac12-\alpha$.
\item[(ii)] If for some $\alpha\in(0,1/2)$,
%
\begin{equation}
\label{ordre4} \frac{d^*\log(d/d^*)}{n}\geq\frac{\kappa}{\alpha},
\end{equation}
then $\inf_{\widetilde{J}}\sup_{\sff\in\Sigma} \Pb_\sff
(\widetilde
{J}\neq J_\sff )\geq\frac12-\alpha$.
\end{longlist}
\end{theorem}

\begin{figure}

\includegraphics{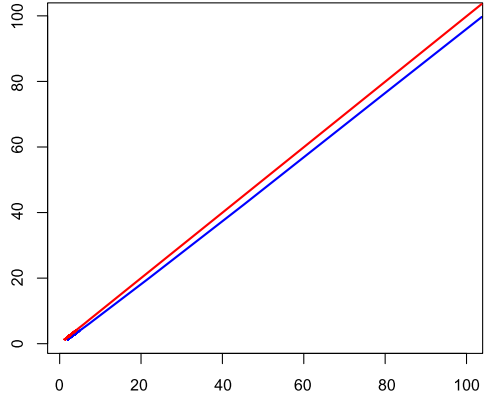}

\caption{The curve of the function $L\mapsto{\gammaLL}$ (blue) and the
bisector (red).}\label{fig2}
\end{figure}

It is worth stressing here that condition (\ref{ordre3}) is the
converse of condition (\ref{ordre2}) of
Theorem~\ref{prop2} in the case $d^*\to\infty$, in the sense that
condition (\ref{ordre2}) amounts to
requiring that the left-hand side of (\ref{ordre3}) is smaller than
some constant. There is, however, one
difference between the quantities involved in these conditions: the
term $N(d^*,\vartheta(1+\tau))$ of (\ref{ordre2})
is replaced by $N(d^*,{\gammaL})$ in condition (\ref{ordre3}). One can
wonder how close
${\gammaL}$ is to $\vartheta$. To give a qualitative answer to this
question, we plotted in Figure~\ref{fig2}
the curve of the mapping $\vartheta\mapsto{\gammaL}$ along with the
bisector $\vartheta\mapsto\vartheta$.
We observe that the difference between two curves is small compared to
$\vartheta$. As we discuss it later, this property shows that
the constants involved in the necessary condition and in the sufficient
condition for consistent estimation
of $J$ are very close, especially for large values of $\vartheta$.

\section{Adaptivity and minimax rates of separation}\label{secadapt}

\subsection{\texorpdfstring{Adaptation with respect to $L$ and $\kappa$}{Adaptation with respect to L and kappa}}\label{ssecadapt}
The estimator $\widehat J(\mm,\llambda)$ we have introduced in
Section~\ref{sec3} is clearly nonadaptive:
the tuning parameters $(\mm,\llambda)$ recommended by the developed
theory involve the values $L$ and $\kappa$, which
are generally unknown.\vadjust{\goodbreak} Fortunately, we can take advantage of the fact
that the choice of $\mm$ and~$\llambda$ is
governed by the one-dimensional parameter $\vartheta=L(1+\tau)/\kappa$.
Therefore, it is realistic to assume
that a finite grid of values $1<\vartheta_1\le\cdots\le\vartheta_{K}<\infty$ is available containing a true
value of $\vartheta$. The following result provides an adaptive
procedure of variable selection with
guaranteed control of the error.

\begin{proposition}\label{propadapt}
Let $1<\vartheta_1\le\cdots\le\vartheta_{K}<\infty$ and $\tau>0$ be
given values, and se\setcounter{footnote}{1}t\footnote{We use the convention
that the minimum over an empty set equals $+\infty$.}
\[
i^*=\min \biggl\{i\dvtx (1+\tau)\frac{\max_{j=1,\ldots,d}
\sum_{\bk\in\ZZ^d} k_j^2\theta_\bk^2}{\min_{j\in J}\sum_{\bk
\dvtx k_j\neq0}
\theta_\bk^2}\le\vartheta_i
\biggr\}\le K.
\]
For every $i,\ell\in\NN$, let us denote $\widehat J_n(i)=\widehat
J_n(\mm(\vartheta_i),\llambda(\vartheta_i))$
with $m_\ell(\vartheta)=(\vartheta\ell)^{1/2}$ and
\[
\lambda_\ell(\vartheta)=\frac{2\sqrt{2N(\ell,\vartheta)d^*\log
(2ed/d^*)}+4d^*\log(2ed/d^*)}{n}.
\]
If the condition $4\lambda_s(\vartheta_{i^*})<\kappa\tau/(1+\tau)$ is
fulfilled, then the estimator
$\widehat J_n^{\mathrm{ad}}=\bigcup_{i=1}^K \widehat J_n(i)$ satisfies
$\Pb(\widehat J_n^{\mathrm{ad}}\not= J)\le(K+2)(d^*/2ed)^{d^*}$.
\end{proposition}
In simple words, if the grid of possible values $\{\vartheta_i\}$ has a
cardinality $K$ which is not too large
[i.e., $K(d^*/d)^{d^*}\to0$], then declaring a variable relevant if at
least one of the procedures
$\widehat J_n(i)$ suggests its relevance provides a consistent and
adaptive variable selection strategy. The proof
of this statement follows immediately from Proposition~\ref{prop0} and
Theorem~\ref{prop2}. Indeed, applying
Proposition~\ref{prop0} with $A=2$ yields $\Pb(\widehat J_n^{\mathrm{ad}}\not
\subset J)\le
\sum_{i=1}^K \Pb(\widehat J_n(i)\not\subset J)\le K(d^*/2ed)^{d^*}$,
while Theorem~\ref{prop2} ensures that
$\Pb(\widehat J_n^{\mathrm{ad}}\not\supset J)\le\Pb(\widehat
J_n(i^*)\not
\supset J)\le2(d^*/2ed)^{d^*}$.

\subsection{Minimax rates of separation}\label{ssecmmx}
Since the methodology of Section~\ref{sec3} takes its roots in the
theory of hypotheses testing, one
naturally wonders what are the minimax rates of separation in the
problem of variable selection. The results
stated in foregoing sections allow us to answer this question in the
case of Sobolev smoothness 1 and
alternatives separated in $L^2$-norm. The following result, the proof
of which is postponed to the
Appendix~\ref{appprop5} provides minimax rates. We assume herein that
the true sparsity
$s=\operatorname{Card}(J)$ and its known upper estimate $d^*$ are such
that $d^*/s$ is bounded from above by some
constant.

\begin{proposition}\label{propmmx}
$\!\!\!$There is a constant $D^*$ depending only on $L$ such that~if
\[
\kappa\ge D^* \biggl\{ \biggl(\frac{\log(d/s)}{n^2} \biggr)^{2/(4+s)}\vee
\frac{s\log(d/s)}{n} \biggr\},\vadjust{\goodbreak}
\]
then there exists a consistent estimator of $J$. Furthermore, the
consistency is uniform in
$\sff\in\Sigma(\kappa,L)$. On the other hand, there is a constant $D_*$
depending only on $L$
such that if
\[
\kappa\le D_* \biggl\{ \biggl(\frac{\log(d/s)}{n^2} \biggr)^{2/(4+s)}\vee
\frac{s\log(d/s)}{n} \biggr\},
\]
then uniformly consistent estimation of $J$ is impossible.
\end{proposition}

Borrowing the terminology of the theory of hypotheses testing, we say that
$ (\frac{\log(d/s)}{n^2} )^{2/(4+s)}\vee\frac{s\log(d/s)}{n}$
is the minimax rate of separation in
the problem of variable selection for Sobolev smoothness one. These
results readily extend to Sobolev smoothness
of any order $\beta\ge1$, in which case the rate of separation takes
the form
$ (\frac{\log(d/s)}{n^2} )^{2\beta/(4\beta+s)} \vee\frac
{s\log
(d/s)}{n}$.
The first term in this maximum coincides, up to the logarithmic term,
with the minimax rate of separation
in the problem of detection of an $s$-dimensional signal \cite
{IngsterStepanova11}. Note, however, that in our case this logarithmic
inflation
is unavoidable. It is the price to pay for not knowing in advance
which $s$ variables are
relevant.

\section{Nonparametric regression with random design}\label{sec6}

So far, we have analyzed the situation in which noisy observations of
the regression function $\sff(\cdot)$ are
available at all points $\mathbf{x}\in[0,1]^d$. Let us turn now to
the more
realistic model of nonparametric regression,
when the observed noisy values of $\sff$ are sampled at random in the
unit hypercube $[0,1]^d$.
More precisely, we assume that $n$ independent and identically
distributed pairs of input-output variables $(\bX_i,Y_i)$,
$i=1,\ldots,n$ are observed that obey the regression model
\[
Y_i=\sff(\bX_i)+\sigma\varepsilon_i, \qquad i=1,
\ldots,n.
\]
The input variables $\bX_1,\ldots,\bX_n$ are assumed to take values in
$\RR^d$ while the output variables $Y_1,\ldots,Y_n$
are scalar. As usual, $\varepsilon_1,\ldots,\varepsilon_n$ are such
that $\Ex[\varepsilon_i|\bX_i]=0$, $i=1,\ldots,n$; additional conditions
will be imposed later. Without requiring from $\sff$ to be of a special
parametric form, we aim at recovering
the set $J\subset\{1,\ldots,d\}$ of its relevant variables. The noise
magnitude $\sigma$ is assumed to be
known.

It is clear that the estimation of $J$ cannot be accomplished without
imposing some further assumptions on $\sff$ and on the
distribution $P_X$ of the input variables. Roughly speaking, we will
assume that $\sff$ is differentiable with a squared
integrable gradient and that $P_X$ admits a density which is bounded
from below. More precisely, let $\sfg$ denote the
density of $P_X$ w.r.t. the Lebesgue measure.
\begin{longlist}[{[C2]}]
\item[{[C2]}] $\sfg(\mathbf{x})=0$ for any $\mathbf{x}\notin
[0,1]^d$ and that $\sfg
(\mathbf{x})\ge{\gmin}$ for any $\mathbf{x}\in[0,1]^d$.
\end{longlist}

The next assumptions imposed to the regression function and to the
noise require their boundedness in an appropriate sense.\vadjust{\goodbreak}
These assumptions are needed in order to prove, by means of a
concentration inequality, the closeness of the empirical
coefficients to the true ones.
\begin{longlist}[{[C3]}]
\item[{[C3$(L_\infty,L_2)$]}] The $L^\infty([0,1]^d,\RR,P_X)$ and
$L^2([0,1]^d,\RR,P_X)$ norms of the function $\sff$ are
bounded from above, respectively, by $L_\infty$ and $L_2$, that is,
$\Pb (|\sff(\bX)|\le L_\infty )=1$ and $\Ex[\sff(\bX
)^2]\le L_2^2$.
\item[{[C4]}] The noise variables satisfy a.e. $\Ex[e^{t \varepsilon
_i}|\bX_i]\le e^{t^2/2}$ for all $t>0$.
\end{longlist}

We stress once again that the primary aim of this work is merely to
understand when it is possible to consistently estimate the
sparsity pattern. The estimator that we will define is intended to show
the possibility of consistent estimation, rather than
being a practical procedure for recovering the sparsity pattern.
Therefore, the estimator will be allowed to depend on the
parameters $\gmin$, $L$, $\kappa$ and $L_2$ appearing in conditions
[C1]--[C3].

\subsection{An estimator of $J$ and its consistency}
The estimator of the sparsity pattern $J$ that we are going to
introduce now is based on the following simple observation:
if $j\notin J$, then $\theta_\bk[\sff]=0$ for every $\bk$ such that
$k_j\not=0$. In contrast, if $j\in J$, then there exists
$\bk\in\ZZ^d$ with $k_j\not=0$ such that $|\theta_\bk[\sff]|>0$. To
turn this observation into an estimator of $J$, we
start by estimating the Fourier coefficients $\theta_\bk[\sff]$ by
their empirical counterparts,
\[
\hatthetak=\frac{1}{n}\sum_{i=1}^n
\frac{\varphi_\bk(\bX_i)}{\sfg(\bX_i)}Y_i,\qquad \bk\in\ZZ^d.
\]
Then, for every $\ell\in\NN$ and for any $\gamma>0$, we introduce the
notation $S^j_{m,\ell}=
\{\bk\in\ZZ^d\dvtx \|\bk\|_{2}\le m, \|\bk\|_0\le\ell,
k_j\neq0
\}$. The estimator of $J$ is
defined by
%
\begin{equation}
\label{hatJ} \widehat{J}_n^{(1)}(m,\lambda)= \bigl\{j\in
\{1,\ldots,d\}\dvtx \operatorname {max}_{\bk\in S_{m,d^*}^j} |\hatthetak|>{\lambda}
\bigr\},
\end{equation}
where $m$ and ${\lambda}$ are some parameters to be defined later. The
next result, the proof of which is placed in the
supplementary material~\cite{CDsupp}, provides consistency guarantees
for $\widehat{J}_n^{(1)}(m,\lambda)$.

\begin{theorem}\label{thm1}
Let conditions \textup{[C1]--[C4]} be fulfilled with some known values
$\gmin
$, $\vartheta=2L/\kappa$ and $L_2$.
Assume furthermore that the design density $\sfg$ and an upper estimate
on the noise magnitude $\sigma$ are
available. Set $m=(\vartheta d^*)^{1/2}$ and $\lambda=4(\sigma
+L_2)
({d^*\log(24 \sqrt{\vartheta}d/d^*)}/{n\gmin^2})^{1/2}$.
If the following conditions are satisfied:
%
\begin{eqnarray}\label{cond3}
\frac{d^*\log(24\sqrt{\vartheta} d/d^*)}{n}&\leq&\frac
{L_2^2}{L_\infty^2},
\nonumber
\\[-8pt]
\\[-8pt]
\nonumber
 \frac{128({\sigma}+L_2)^2d^*N(d^*,\vartheta)\log(24\sqrt{\vartheta}
d/d^*)}{n\gmin^2}&<& \kappa,
\end{eqnarray}
then the estimator $\widehat{J}^{(1)}(m,\lambda)$ satisfies $\Pb
(\widehat{J}^{(1)}(m,\lambda)\neq J ) \leq(8d/d^*)^{-d^*}$.\vadjust{\goodbreak}
\end{theorem}

If we take a look at the conditions of Theorem~\ref{thm1} ensuring the
consistency of $\widehat J_n^{(1)}$,
it becomes clear that the strongest requirement is the second
inequality in~(\ref{cond3}). Roughly speaking,
this condition requires that $d^*N(d^*,\vartheta)\log(d/d^*)/n$ is
bounded from above by some constant.
According to results stated in Section~\ref{sec4}, $N(d^*,\vartheta)$
diverges exponentially fast,
making inequality (\ref{cond3}) impossible for $d^*$ larger than $\log
n$ up to a multiplicative constant.

It is also worth stressing that although we require the $P_X$-a.e. boundedness of~$\sff$ by some constant $L_\infty$,
this constant is not needed for computing the estimator proposed in
Theorem~\ref{thm1}. Only constants related to some quadratic
functionals of the sequence of Fourier coefficients $\theta_\bk[\sff]$
are involved in the tuning parameters $m$ and $\lambda$.
This point might be important for designing practical estimators of
$J$, since the estimation of quadratic functionals is more
realistic (see, e.g.,~\cite{LaurentMassart,Cai06}) than the
estimation of $\sup$-norm.

Theorem~\ref{thm1} can be reformulated to characterize the level of
relevance $\kappa$ for the relevant components of $\bX$
making their identification possible. In fact, an alternative way of
stating Theorem~\ref{thm1} is the following:
under conditions [C1]--[C4] if $\vartheta$ is an arbitrary
tuning parameter satisfying the first inequality
in (\ref{cond3}), then the estimator $\widehat J_n^{(1)}(m,\lambda
)$---with $m$ and $\lambda$ chosen as in
Theorem~\ref{thm1}---satisfies $\Pb(\widehat J_n^{(1)}(m,\lambda
)\not
=J)\le(8d/d^*)^{-d^*}$ if the smallest
level of relevance $\kappa$ for components $X_j$ of $\bX$ with $j\in J$
is not smaller than $8\lambda^2 N(d^*,m^2/d^*)$.
This statement can be easily deduced from the proof of Theorem \ref
{thm1}; cf. the supplementary material~\cite{CDsupp}.\vspace*{-2pt}

\subsection{Tightness of the assumptions}
A natural question is now to check that the assumptions of Theorem \ref
{thm1} are tight in the asymptotic
regimes of fixed sparsity and increasing ambient dimension, as well as
increasing sparsity.
We will only establish an analogue of claim (ii) of Theorem \ref
{thm3}. An attempt to prove a result similar
to claim (i) of Theorem~\ref{thm3} was done in~\cite{ComDal11}, Theorem 2. However, the result
of~\cite{ComDal11} involves a stringent assumption on the empirical
Gram matrix (cf. condition (6) in~\cite{ComDal11})
and, unfortunately, we are unable to prove the existence of a sampling
scheme for which this assumption is fulfilled.

We assume that the errors $\varepsilon_i$ are i.i.d. standard
Gaussian, and
we focus our attention on the functional class $\Sigma(\kappa,L)$. The
following simple result shows that
the conditions of Theorem~\ref{thm1} are tight in the case of fixed
intrinsic dimension.\vspace*{-2pt}

\begin{proposition}\label{prop4}
Let the design $\bX_1,\ldots, \bX_n\in[0,1]^d$ be either deterministic
or random. If for some
positive $\alpha<1/2$, the inequality
\[
\frac{d^*\log(d/d^*)}{n}\geq\kappa\alpha^{-1}
\]
{\spaceskip=0.19em plus 0.05em minus 0.02em holds true, then there is a constant $c\!>\!0$ such that \mbox{$\inf_{\widetilde J_n}\sup_{\sff\in\Sigma(\kappa,L)}
\!\Pb_\sff(\widetilde J_n\!\neq\! J_\sff)\!\geq\! c$}.}\vadjust{\goodbreak}
\end{proposition}

\section{Concluding remarks}\label{sec7}

The results proved in previous sections almost exhaustively answer the
questions on the existence of consistent
estimators of the sparsity pattern in the model of Gaussian white noise
and, to a smaller extent, in nonparametric
regression. In fact as far as only rates of convergence are of
interest, the result obtained in Theorem~\ref{prop2}
is shown in Section~\ref{sec5} to be unimprovable. Thus only the
problem of finding sharp constants remains open.
To make these statements more precise, let us consider the simplified
set-up $\sigma=\kappa=1$ and define the
following two regimes:
\begin{itemize}
\item The regime of fixed sparsity, that is, when the sample
size $n$ and the ambient dimension $d$
tend to infinity but the intrinsic dimension $d^*$ remains constant or bounded.
\item The regime of increasing sparsity, that is, when the
intrinsic dimension $d^*$ tends to infinity
along with the sample size $n$ and the ambient dimension $d$. For
simplicity, we will assume that $d^*=O(d^{1-\varepsilon})$ for some
$\varepsilon>0$.
\end{itemize}
In the fixed sparsity regime, in view of Theorems~\ref{prop2} and \ref
{thm1}, consistent estimation of the sparsity pattern
can be achieved both in the Gaussian white noise model and
nonparametric regression as soon as $\limsup_{n\to\infty}(d^*\log d)/n
< c_\star$, where
$c_\star$ is the constant defined by $c_*=1/8$ for the Gaussian white
noise model and
\[
c_\star=\min \biggl(\frac{L_2^2}{2L_\infty^2}, \frac
{\gmin^2}{2^8(1+L_2)^2N(d^*,2L)} \biggr)
\]
for the regression model. On the other hand, by Theorem~\ref{thm3}
and Proposition~\ref{prop4}, consistent estimation of the sparsity
pattern is impossible if\break
$\liminf_{n\to\infty}(d^*\log d)/n> c^\star$ with $c^\star=2$.
Thus, up
to multiplicative constants
$c_\star$ and $c^\star$ (which are clearly not sharp), the results of
Theorems~\ref{prop2} and~\ref{thm1}
cannot be improved in the regime of fixed sparsity.

\begin{figure}

\includegraphics{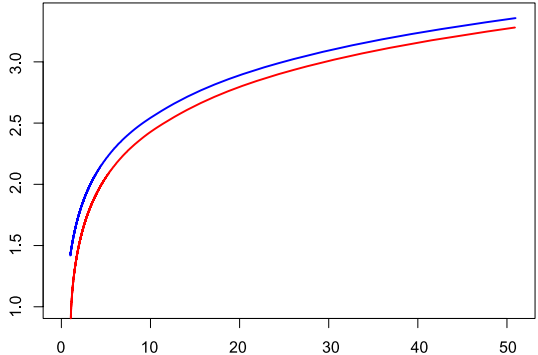}

\caption{The curves of functions $L\mapsto\sfl_{L}(z_L)$ (blue curve)
and $L\mapsto\sfl_{{\gammaLL}}(z_{{\gammaLL}})$ (red curve).}\label{fig4}
\end{figure}

In the regime of increasing sparsity, the results we get in the model
of Gaussian white noise are much stronger than those for nonparametric
regression. In the former model, taking the logarithm of both sides of
inequality (\ref{ordre2}) and using formula~(\ref{eqlog}) for
$N(d^*,\cdot)=N_1(d^*,\cdot)-N_2(d^*,\cdot)$, we see that consistent
estimation of $J$ is possible when, for some $\tau>0$ and for all $n$,
the following two conditions are fulfilled:
%
\begin{eqnarray}
\label{eq10} %
\cases{ \sfl_{L+\tau}(z_{L+\tau}) d^*+
\frac12\log d^*+\log\log \bigl(d/d^*\bigr)-2\log n<\underline{c}_1,
\vspace*{2pt}
\cr
\log d^*+\log\log\bigl(d/d^*\bigr)-\log n \le
\underline{c}'_1 } %
\end{eqnarray}
with some constants $\underline{c}_1=\underline{c}_1(L,\tau)$ and
$\underline{c}_1'=\underline{c}_1'(L,\tau)$.
On the other hand, Theorem~\ref{thm3} yields that there are some
constants $\bar c_1$ and
$\bar c_1'$ such that it is impossible to consistently estimate $J$ if
either one of the conditions
%
\begin{eqnarray}
\label{eq11} \sfl_{{\gammaLL}}(z_{{\gammaLL}}) d^*+\tfrac12\log d^*+\log\log
\bigl(d/d^*\bigr)-2\log n&\ge&\bar{c}_1,
\\
\label{eq12b} \log d^*+\log\log\bigl(d/d^*\bigr)-\log n &\ge&
\bar{c}'_1,
\end{eqnarray}
is satisfied.
First note that the left-hand side of the second condition in (\ref
{eq10}) is exactly the same as the left-hand side of (\ref{eq12b}).
If we compare now the left-hand side of the first condition in (\ref
{eq10}) with the left-hand side of (\ref{eq11}), we see that
only the coefficients of $d^*$ differ. To measure the degree of
difference of these two coefficients we draw in Figure~\ref{fig4}
the plots of the functions $L\mapsto\sfl_L(z_{L})$ and $L\mapsto\sfl_{{\gammaLL}}(z_{{\gammaLL}})$, with ${\gammaLL}$ as
is Theorem~\ref{thm3}. One can observe that the two curves are very
close, especially for relatively large values of $L$. This implies
that the conditions in (\ref{eq10}) are tight.
A~simple consequence of inequalities (\ref{eq10}) and (\ref{eq11})
is that
the consistent recovery of the sparsity pattern is possible under the
condition $d^*/\log n\to0$ and impossible for $d^*/\log n\to\infty$ as
$n\to\infty$, provided that $\log\log(d/d^*)=o(\log n)$.

Still in the regime of increasing sparsity, but for nonparametric
regression, we proved that consistent estimation
of the sparsity pattern is possible whenever
%
\begin{eqnarray}
\label{eq10b} %
\cases{ \sfl_{L+\tau}(z_{L+\tau}) d^*+
\frac12\log d^*+\log\log \bigl(d/d^*\bigr)-\log n<\underline{c}_2,
\vspace*{2pt}
\cr
\log d^*+\log\log d-\log n < \underline{c}'_2
} %
\end{eqnarray}
with some constants $\underline{c}_2=\underline{c}_2(\gmin,\sigma
,L_2,L)$ and $\underline{c}'_2=2\log(L_2/L_\infty)$.
As we have already mentioned, the second condition in (\ref{eq10b})
is tight, up to the choice of $\underline{c}_2'$, in view of
Proposition~\ref{prop4}. It is natural to expect that the first
condition is
tight as well, since it is in the model of Gaussian white noise, which
has the reputation of being simpler than the model of nonparametric
regression. However, we do not have a mathematical proof of this statement.

Let us stress now that, all over this work, we have deliberately
avoided any discussion on the computational aspects of the variable selection
in nonparametric regression. The goal in this paper was to investigate
the possibility of consistent recovery without paying attention to the
complexity of the selection procedure. This lead to some conditions
that could be considered a benchmark for assessing the properties of
sparsity pattern estimators. As for the estimators proposed in
Section~\ref{sec3}, it is worth noting that their computational
complexity is not
always prohibitively large. A recommended strategy is to compute the
coefficients $\hatthetak$ in a stepwise manner; at each step
$K=1,2,\ldots,d^*$ only the coefficients $\hatthetak$ with $\|\bk\|_0=K$ need to be computed and compared with the threshold. If some
$\hatthetak$ exceeds the threshold, then all the variables $X^j$
corresponding to nonzero coordinates of $\bk$ are considered as relevant.
We can stop this computation as soon as the number of variables
classified as relevant attains $d^*$. While the worst-case complexity
of this
procedure is exponential, there are many functions $\sff$ for which the
complexity of the procedure will be polynomial in $d$. For example, this
is the case for additive models in which $\sff(\mathbf{x})= \sff_1(x_{i_1})+\cdots+\sff_{d^*}(x_{i_{d^*}})$ for some univariate functions
$\sff_1,\ldots,\sff_{d^*}$.

Note also that in the present study we focused exclusively on the
consistency of variable selection without paying any attention to
the consistency of regression function estimation. A thorough analysis
of the latter problem being left to a future work, let us simply remark
that in the case of fixed $d^*$, under the conditions of Theorem \ref
{thm1}, it is straightforward to construct a consistent estimator of
the regression function. In fact, it suffices to use a projection
estimator with a properly chosen truncation parameter on the set of relevant
variables. The situation is much more delicate in the case when the
sparsity $d^*$ grows to infinity along with the sample size $n$. Presumably,
condition (\ref{eq10}) is no longer sufficient for consistently
estimating the regression function. The rationale behind this
conjecture is that
the minimax rate of convergence for estimating $\sff$ in our context,
if we assume in addition that the set of relevant variables is known, is
equal to $n^{-2/(2+d^*)}=\exp(-2\log n/(2+d^*))$. If the left-hand side
of (\ref{eq10}) is equal to a constant and $\log\log d=o(\log n)$, then
the aforementioned minimax rate does not tend to zero, making thus the
estimator inconsistent.

Finally, we would like to mention that the selection of relevant
variables is a challenging statistical task, which might be useful to
perform independently of the task of regression function estimation.
Indeed, if we succeed in identifying relevant variables on a data-set
having a small sample size, we can continue the data collection process
more efficiently by recording only the values of relevant variables.
This may considerably reduce the memory costs related to the data
storage and the financial costs necessary for collecting new data. Then,
the regression function may be estimated more accurately on the base of
this new (larger) data-set.

\begin{appendix}\label{app}

\section{\texorpdfstring{Proof of Proposition \lowercase{\protect\ref{prop0}}}
{Proof of Proposition 1}}\label{appA0}

To ease notation, we write $\widehat J_n$ instead of $\widehat
J_{n}(\mm
,\llambda)$. It is clear that
$\widehat J_n\not\subset J$ if and only if $\exists j\in J^c$ such that
$\max_{\ell\le d^*}\lambda_\ell^{-1}\max_{I\in P_\ell^d} \widehat
Q_{m,I}^j\ge1$,
where $Q_{m,I}^j=\sum_{\bk\in S_{m,I}^j} \theta_\bk^2$. For every\vadjust{\goodbreak}
$j\in
\{1,\ldots,d\}$, let us set
$R_{m,I}^j=\sum_{\bk\in S_{m,I}^j}(\xi_\bk^2-1)$ and
$N_{m,I}^j=(Q_{m,I}^j)^{-1/2}\sum_{\bk\in S_{m,I}^j} \theta_\bk\xi_\bk$
so that
%
\begin{equation}
\label{notation1} \widehat{Q}_{m,I}^j=\sum
_{\bk\in S_{m,I}^j} \biggl(y_\bk^2-
\frac
{1}{n} \biggr)= Q_{m,I}^j+\frac{2\sqrt{Q_{m,I}^j}}{\sqrt{n}}
N_{m,I}^j+\frac1n R_{m,I}^j.
\end{equation}
For $j\in J^c$, the first two terms of the last sum vanish and,
therefore, we have
\[
\{\widehat J_n\not\subset J \}= \bigcup
_{j\in J^c}\bigcup_{\ell
\le d^*}\bigcup
_{I\in P_\ell^d} \bigl\{R_{m,I}^j\ge n
\lambda_\ell \bigr\}=\bigcup_{\ell\le
d^*}\bigcup
_{I\in P_\ell^d}\bigcup_{j\in J^c\cap I}
\bigl\{R_{m,I}^j\ge n\lambda_\ell \bigr\},
\]
where the last equality\vspace*{-2pt} results from the fact that $R^j_{m,I}=0$ if
$j\notin I$. The random variable $R^j_{m,I}$,
being a centered sum of squares of independent standard Gaussian random
variables, follows a translated
$\chi^2$-distribution. The tails of this distribution can be evaluated
using the following result.
%
\begin{lemma}[(cf. Lemma 1 in~\cite{LaurentMassart})]\label{lemmassart}
Let $\xi_1,\ldots,\xi_D$ be independent standard Gaussian random variables.
For every $x\ge0$ and for every vector $\ba=(a_1,\ldots,a_D)\in\RR^D_+$,
the following inequalities
hold true:
\begin{eqnarray*}
\Pb \Biggl(\sum_{i=1}^D
a_i\bigl(\xi_i^2-1\bigr)\geq2\|\ba
\|_2\sqrt{x}+2\|\ba \|_\infty x \Biggr)&\leq&\exp(-x),
\\
\Pb \Biggl(\sum_{i=1}^D
a_i\bigl(\xi_i^2-1\bigr)\leq-2\|\ba
\|_2\sqrt{x} \Biggr)&\leq &\exp(-x).
\end{eqnarray*}
\end{lemma}
We apply this lemma to $R^i_{m_\ell,I}$, for which $\|\ba\|_\infty=1$
and $\|\ba\|_2^2={N(\ell,m_\ell^2/\ell)}$.
Setting $n\lambda_\ell=2 \sqrt{N(\ell,m_\ell^2/\ell)x}+2x$ and using
the union bound, we get
\begin{eqnarray*}
\Pb (\widehat J_n\not\subset J )&\le& \Pb \Biggl(\bigcup
_{\ell=1}^{d^*} \Bigl\{\max_{I\in P_\ell^d;i\in I}
R^i_{m_\ell,I} \ge n\lambda_\ell \Bigr\} \Biggr)
\\
&\leq&\sum_{\ell=1}^{d^*} \ell
\operatorname{Card}\bigl(P_\ell^d\bigr) \max_{I\in P_\ell^d;i\in I}
\Pb \bigl( R^i_{m_\ell,I} >n\lambda_\ell \bigr) \leq
e^{-x}\sum_{\ell=1}^{d^*}\ell
\pmatrix{{d}
\cr
{\ell }}.
\end{eqnarray*}
One checks that $\sum_{\ell=1}^{d^*} \ell({{d}\atop{\ell}}) \le
(2ed/d^* )^{d^*}$ holds true
for every pair of integers $(d^*,d)$ such that $1\le d^*\le d$; cf. the supplementary material~\cite{CDsupp} for a proof.
Hence, for $x=Ad^*\log(2ed/d^*)$, we get $\Pb (\widehat J_n\not
\subset J )\le(2ed/d^*)^{-(A-1)d^*}$.

\section{\texorpdfstring{Proof of Theorem \lowercase{\protect\ref{prop2}}}
{Proof of Theorem 1}}\label{appA1}

We begin with proving a stronger result that implies the claim of
Theorem~\ref{prop2}.

\begin{proposition}\label{propnew}
Let $\alpha$ be a real number from $(0,1)$. If for every $j\in J$ and
for $s=\operatorname{Card}(J)$ the inequality
%
\begin{equation}
\label{QQ1}  Q_{m_s,J}^j \ge \biggl\{ \biggl[
\lambda_s+\frac{2\sqrt {N(s,m_s^2/s)\log
(2s/\alpha)}+1}{n} \biggr]^{1/2}+ \biggl[
\frac{2\log(2s/\alpha)}{n} \biggr]^{1/2} \biggr\}^2\hspace*{-35pt}
\end{equation}
holds true, then $\Pb(J\not\subset\widehat J_n)\le\alpha$.
\end{proposition}
\begin{pf}
To bound from above the probability of type II error, we rely on~the~equivalence:
$J\not\subset\widehat J_n$ if and only if $\exists j\in J$ such
that\break
$\max_{\ell\le d^*}\lambda_\ell^{-1}\max_{I\in P_\ell^d} \widehat
Q_{m,I}^j\le1$.
Recall that $s=\operatorname{Card}(J)$. Using Bonferroni's inequality,
we get
%
\begin{eqnarray}
\label{bonf2} \Pb (J\not\subset\widehat J_n ) &\le&\sum
_{j\in J} \Pb \Bigl(\max_{\ell\le d^*}\lambda_\ell^{-1}
\max_{I\in P_\ell^d} \widehat Q_{m,I}^j\le1 \Bigr)
\nonumber
\\[-8pt]
\\[-8pt]
\nonumber
&\le&\sum_{j\in J} \Pb \bigl(\widehat
Q_{m_{s},J}^j\le\lambda_s \bigr) \le s
\max_{j\in J}\Pb \bigl(\widehat Q_{m_{s},J}^j\le
\lambda_s \bigr).
\end{eqnarray}
By virtue of decomposition (\ref{notation1}),
\[
\Pb \bigl(\widehat Q_{m_s,J}^j\le
\lambda_s \bigr) 
=\Pb \biggl( \biggl({ \sqrt
{Q_{m_s,J}^j}}+\frac{1}{\sqrt {n}}N_{m_s,J}^j
\biggr)^2+\frac1n \bigl(R_{m_s,J}^j-
\bigl(N_{m_s,J}^j\bigr)^2 \bigr)\le
\lambda_s \biggr).
\]
One checks that $R_{m_s,J}^j-(N_{m_s,J}^j)^2+N(s,m_s^2/s)$ is a drawn
from $\chi^2$-distribu\-tion~with
$N(s,m_s^2/s)-1$ degrees of freedom. Therefore, using Lemma \ref
{lemmassart} stated in previous section, we get
$
\Pb (\frac1n  (R_{m_s,J}^j-(N_{m_s,J}^j)^2 )+\frac1n \le\break
-2\sqrt{N(s,m_s^2/s)\log(2s/\alpha)}  )
\le\frac{\alpha}{2s}$.
Therefore, $\Pb (\widehat Q_{m_s,J}^j\le\lambda_s )$ is
upper-bounded by
\begin{eqnarray*}
\frac{\alpha}{2s}+ \Pb \biggl( \biggl({ \sqrt{Q_{m_s,J}^j}}+
\frac
{1}{\sqrt {n}}N_{m_s,J}^j \biggr)^2 \le
\lambda_s+\frac{2\sqrt{N(s,m_s^2/s)\log(2s/\alpha)}+1}n \biggr).
\end{eqnarray*}
Using the condition of the proposition, we get
$\Pb (\widehat Q_{m_s,J}^j\le\lambda_s )
\le\frac{\alpha}{2s}+\break \Pb ( N_{m_s,J}^j\le-\sqrt{2\log
(2s/\alpha
)} )\le\frac{\alpha}{s}$.
Combining this inequality with (\ref{bonf2}), we get the result of
Proposition~\ref{propnew}.
\end{pf}

To deduce the claim of Theorem~\ref{prop2} from that of
Proposition~\ref{propnew}, we use the following
lower bound:
%
\begin{eqnarray}
\label{minor} Q_{m_s,J}^{j}&=&Q^{j}-\sum
_{j\in\supp(\bk)\subset J} \theta_{\bk
}^2\1_{\{
\|\bk\|_2\ge m_s\}}
\geq\kappa-\sum_{j\in\supp(\bk)\subset J} \theta_{\bk}^2
\1_{\{\|
\bk\|
_2\ge m_s\}}
\nonumber
\\[-8pt]
\\[-8pt]
\nonumber
&\geq&\kappa- m_s^{-2}\sum_{j\in\supp(\bk)\subset J}
\theta_{\bk
}^2\| \bk\|_2^2 \geq
\kappa- m_s^{-2} Ls
\end{eqnarray}
for every $j\in J$. Our choice of $m_s$, $m_s=\sqrt{sL(1+\tau)/\kappa
}$, ensures that $Q_{m_s,J}^{j}\ge\kappa\tau/(1+\tau)$.
Finally, using a very rough bound (which is sufficient for our
purposes), the right-hand side in (\ref{QQ1}) can be
upper-bounded by $4\lambda_s$ if $\alpha$ is chosen to be equal to
$2(2ed/d^*)^{-(A-1)d^*}$. Therefore, if
$\frac{\kappa\tau}{1+\tau} \ge4\lambda_s,$ then (\ref{QQ1}) holds
true with $\alpha=2(2ed/d^*)^{-(A-1)d^*}$ and, therefore,
the type II error has a probability less than or equal to
$2(2ed/d^*)^{-(A-1)d^*}$.

\section{\texorpdfstring{Proof of Proposition \lowercase{\protect\ref{prop1}}}
{Proof of Proposition 2}}
\textit{Proof of the first assertion}. This proof can be found
in~\cite{Mazo}; we repeat here the arguments therein for the sake of keeping
the paper self-contained. Recall that $N_1(d^*,\gamma)$ admits an
integral representation with the integrand
\[
\frac{\sfh(z)^{d^*}}{z^{\gamma d^*}} \frac{1}{z(1-z)}=\frac
{1}{z(1-z)}\exp \biggl[d^*\log
\biggl(\frac{\sfh(z)}{z^\gamma} \biggr) \biggr].
\]
For any $y>0$, we define $\phi(y)=e^{-y}\sfh'(e^{-y})/\sfh
(e^{-y})=\sum_{k\in\ZZ} k^2 e^{-yk^{2}}/ \break\sum_{k\in\ZZ} e^{-yk^{2}}$ in such a way
that
\[
\phi(y)=\gamma\quad\Longleftrightarrow\quad\frac{\sfh'(e^{-y})}{\sfh
(e^{-y})}=\frac{\gamma}{e^{-y}}
\quad\Longleftrightarrow\quad\sfl_\gamma'\bigl(e^{-y}
\bigr)=0.
\]
By virtue of the Cauchy--Schwarz inequality, it holds that $
\sum\! k^{4}e^{-yk^{2}}\sum\! e^{-yk^{2}}> (\sum k^{2}e^{-yk^{2}}
)^2$, $\forall y\in(0,\infty)$,
implying that $\phi'(y)<0$ for all $y\in(0,\infty)$, that is, $\phi$
is strictly decreasing. Furthermore, $\phi$ is
obviously continuous with $\lim_{y\to0}\phi(y)=+\infty$ and $\lim_{y\to
\infty} \phi(y)=0$. These properties imply the
existence and the uniqueness of $y_\gamma\in(0,\infty)$ such that
${\phi
(y_\gamma)}=\gamma$. Furthermore, as the inverse
of a decreasing function, the function $\gamma\mapsto y_\gamma$ is
decreasing as well. We set $z_\gamma=e^{-y_{\gamma}}$
so that $\gamma\mapsto z_\gamma$ is increasing.
We also have
\begin{eqnarray*}
\sfl_\gamma''(\zb) &=& \frac{\sfh''\sfh-(\sfh')^2}{\sfh^2}(\zb
)+\frac
{\gamma}{\zb^2}=z_\gamma^{-2} \biggl\{ \frac{\sum_k (k^4-k^2)\zb^{k^2}}{\sum_k\zb^{k^2}}-
\biggl(\frac
{\sum_k
k^2\zb^{k^2}}{\sum_k\zb^{k^2}} \biggr)^2+\gamma \biggr\}
\\
&=&\zb^{-2} \bigl\{-\phi'(y_\gamma)-
\phi(y_\gamma)+\gamma \bigr\} =-\zb^{-2}\phi'(y_\gamma)>0.
\end{eqnarray*}

\textit{Proof of the second assertion}.
We apply the saddle-point method to the integral representing $N_1$;
see, for example, Chapter IX in~\cite{Dieudonne}.
It holds that
%
\begin{eqnarray}
N_1\bigl(d^*,\gamma\bigr) &=& \frac{1}{2\pi i}\oint_{|z|=\zb}
\frac{\sfh
(z)^{d^*}}{z^{\gamma d^*}} \frac{dz}{z(1-z)}
\nonumber
\\[-8pt]
\\[-8pt]
\nonumber
& =& \frac{1}{2\pi i}\oint_{|z|=\zb}
\bigl\{z(1-z)\bigr\}^{-1}e^{d^*\sfl_\gamma(z)}{\,dz}.
\end{eqnarray}
The first assertion of the proposition provided us with a real number
$\zb$ such that $\sfl_\gamma'(\zb)=0$ and $\sfl_\gamma''(\zb)>0$.
The tangent to the steepest descent curve at $\zb$ is vertical. The
path we choose for integration is the circle with center 0 and
radius $\zb$. As this circle and the steepest descent curve have the
same tangent at $\zb$, applying formula (1.8.1) of
\cite{Dieudonne} [with $\alpha=0$ since $\sfl''(\zb)$ is real and positive],
we get that
\begin{eqnarray*}
&&\oint_{|z|=\zb} \bigl\{z(1-z)\bigr\}^{-1}e^{d^*\sfl_\gamma(z)}{\,dz}\\
&&\qquad=
\sqrt{\frac{2\pi}{d^*\sfl_\gamma''(\zb)}}e^{{{i}\pi/2}}\bigl\{\zb (1-\zb)\bigr
\}^{-1}e^{d^*\sfl_\gamma(\zb)}\bigl(1+o(1)\bigr),
\end{eqnarray*}
when $d^*\to\infty$, as soon as the condition\footnote{$\Re u$ stands
for the real part of the complex number $u$.}
$\Re[\sfl_\gamma(z)-\sfl_\gamma(\zb)]\leq-\mu$
is satisfied for some $\mu>0$ and for any $z$ belonging to the circle
$|z|=|\zb|$ and lying not too close to $\zb$.
To check that this is indeed the case, we remark that $\Re[\sfl_\gamma
(z)]=\log |\frac{\sfh(z)}{z^\gamma} |$.
Hence, if $z=\zb e^{{{i}}\omega}$ with $\omega\in[\omega_0,2\pi
-\omega_0]$ for some $\omega_0\in]0,\pi[$, then
\begin{eqnarray*}
\biggl|\frac{\sfh(z)}{z^\gamma} \biggr|&=&\frac{|1+2z+2\sum_{k> 1}
z^{k^{2}}|}{\zb^\gamma} \le\frac{|1+z|+\zb+2\sum_{k>1}\zb^{k^{2}}}{\zb^\gamma}
\\
&\le&\frac{|1+e^{{{i}}\omega_0}\zb|+\zb+2\sum_{k>1}\zb^{k^{2}}}{\zb^\gamma}.
\end{eqnarray*}
Therefore $\Re[\sfl_\gamma(z)-\Re\sfl_\gamma(\zb)]\leq-\mu$ with
$\mu=\log (\frac{1+2\zb+2\sum_{k > 1}\zb^{k^{2}}}{|1+\zb
e^{{{i}}\omega_0}|+\zb+2\sum_{k > 1}\zb^{k^{2}}} )>0$.
This completes the proof for the term $N_1(d^*,\gamma)$. The term
$N_2(d^*,\gamma)$ can be dealt in the same way.

\section{\texorpdfstring{Proof of Theorem \lowercase{\protect\ref{thm3}}}{Proof of Theorem 2}}
To prove (i) we apply Lemma~\ref{lemfano} with $M=({{d}\atop{d^*}})$ in
conjunction with a standard result, the proof of which can be found
in~\cite{ComDal11} and in the supplementary material~\cite{CDsupp}.

\begin{lemma}\label{lem2}
Let $S$ be a subset of $\ZZ^{d}$ of cardinality ${|S|}$ and $A$
be a constant. Define $\mu_S$ as a discrete measure supported on the
finite set of functions
$\{\sff_{\oomega}= \sum_{\bk\in S}A\omega_\bk\varphi_\bk\dvtx  \oomega
\in\{
\pm1\}^S\}$ such that $\mu_S(\sff=\sff_{\oomega})=2^{-{|S|}}$
for every $\oomega\in\{\pm1\}^S$. 
If we define the probability measure $\PP_S$ by $\PP_S(A)=\int_{\Sigma
(\kappa,L)} \Pb_\sff(A) \mu_S(d\sff)$, for every
measurable set $A\subset\RR^n$, and $\PP_0=\Pb_{\sff_0}$, then
$\KL(\PP_S,\PP_0)\leq|S|A^4 n^2$.
\end{lemma}
Without loss of generality, we can assume $\kappa=1$ (the general case
can be reduced to this one
by replacing $L$ and $n$, respectively, by $L/\kappa$ and $n\kappa$).
Thus, $\vartheta=L$. We denote
the set $\Sigma(1,L)$ by $\Sigma_L$ and choose $\mu_0,\ldots, \mu_M$ as
follows: $\mu_0$ is the Dirac measure
$\delta_0$, $\mu_1$ is defined as in Lemma~\ref{lem2} with
$S=\mathcal
{C}_1(d^*,\gamma_L)$
and $A={ [N(d^*,\gamma_L) ]}^{-1/2}$. The measures $\mu_2,\ldots
,\mu_M$
are defined similarly and correspond to the $M-1$ remaining sparsity
patterns of cardinality $d^*$.

In view of inequality (\ref{minor2}) and Lemma~\ref{lemfano}, it
suffices to show that the measures $\mu_\ell$
satisfy $\mu_\ell(\Sigma_L)=1$ and $\sum_{\ell=0}^M\KL(\PP_\ell
,\PP_0)\le(M+1)\alpha\log M$.
Combining Lemma~\ref{lem2} with $\operatorname
{Card}(S)=N_1(d^*,\gamma_L)$ and inequality (\ref{ordre3}), we get
$
\KL(\PP_\ell,\PP_0)\leq\frac{n^2 N_1(d^*,\gamma_L)}{N(d^*,\gamma_L)^2}
\le\frac{n^2 L}{\gamma_L N(d^*,\gamma_L)}\le\alpha\log M$.
Now, let us show that\break $\mu_1(\Sigma_L)=1$. By symmetry, this will imply
that $\mu_\ell(\Sigma_L)=1$ for
every $\ell$. Since $\mu_1$ is supported by the set $\{\sff_{\oomega}
\dvtx \oomega\in\{\pm1\}^{\mathcal C_1(d^*,\gamma_L)}\}$, it is
clear that $\sum_{k_1\neq0} \theta_\bk^2[\sff_{\oomega}
]=A^2[N_1(d^*,\gamma_L)-N_2(d^*,\gamma_L)]=1$ and
\begin{eqnarray*}
\sum_{\bk\in\ZZ^d} k_j^2
\theta_\bk^2[\sff_{\oomega}]&=&\sum
_{\bk
\in
\mathcal C_1(d^*,\gamma_L)} k_j^2 A^2=
\frac1{d^*} \sum_{j=1}^{d^*}\sum
_{\bk\in\mathcal C_1(d^*,\gamma_L)} k_j^2 A^2\leq
A^2 \gamma_L N_1\bigl(d^*,
\gamma_L\bigr)
\\
&\le&\gamma_L \frac{N_1(d^*,\gamma_L) }{N(d^*,\gamma_L) },\qquad  j=1,\ldots ,d^* .
\end{eqnarray*}
The results stated in Section~\ref{sec4} imply that $N_1(d^*,\gamma_L)/N(d^*,\gamma_L)\sim_{d^*\to\infty} 1+(\sfh(\zb)-1)^{-1}$. Our
choice of $\gamma_L$ ensures that, for $d^*$ large enough, $\sff_{\oomega}
\in\Sigma_L$. This completes the proof of claim (i). To prove
(ii), we still use Lemma~\ref{lemfano} with $\mu_0=\delta_0$ and
$\mu_\ell=\delta_{\sff_{\ell}}$, where
for every $\ell\in\{1,\ldots,M\}$, $\sff_{\ell}$ is chosen as follows.
Let $I_1,\ldots,I_M$ be all
the subsets of $\{1,\ldots,d\}$ containing exactly $d^*$ elements.
We define $\sff_{\ell}$, for $\ell\neq0$, by its Fourier coefficients
$\{
\theta^\ell_\bk\dvtx \bk\in\ZZ^d\}$ as follows:
\[
\theta^\ell_{\bk}= %
\cases{ 1, &\quad $
\bk=(k_1,\ldots,k_d)=(\1_{1\in I_\ell},\ldots,
\1_{d\in
I_\ell}),$\vspace*{2pt}
\cr
0, &\quad  $\mbox{otherwise}.$ } %
\]
Obviously, all the functions $\sff_{\ell}$ belong to $\Sigma$ and,
moreover, each $\sff_\ell$ has $I_\ell$ as sparsity pattern.
One easily checks that our choice of $\sff_\ell$ implies $\KL(\Pb_{\sff
_\ell},\Pb_{\sff_0})=n\|\sff_\ell-\sff_0\|_2^2=n$. Therefore,
if $\alpha\log M =\alpha\log({{d}\atop{d^*}})\ge n$, the desired
inequality is satisfied. To conclude, it suffices
to note that $\log({{d}\atop{d^*}})\ge d^*\log(d/d^*)$.

\section{\texorpdfstring{Proof of Proposition \lowercase{\protect\ref{propmmx}}}
{Proof of Proposition 6}}\label{appprop5}

In view of Theorem~\ref{prop2}, applied with $A=2$ and $\tau=1$, the
consistent [uniformly in $\sff\in\Sigma(\kappa,L)$]
estimation of $J$ is possible if
\[
\frac{8\sqrt{2N(s,2L/\kappa)d^*\log(d/d^*)}+16d^*\log
(d/d^*)}{n}\le
\frac{\kappa}{2}.
\]
Since $d^*/s$ is upper-bounded by some constant, there is a constant
$D^*_1$ such that
the left-hand side of the last display is upper-bounded by
\[
D^*_1 \biggl\{\frac{\sqrt{N(s,2L/\kappa)s\log(d/s)}}{n}\vee\frac
{s\log
(d/s)}{n}
\biggr\}.
\]
As proved in Lemma~\ref{lemcard0} below, $N(s,2L/\kappa)\le0.3
(18\pi e L/\kappa)^{s/2}$. Thus there is a
constant $D_2$ such that
\[
\biggl\{\frac{\sqrt{N(s,2L/\kappa)s\log(d/s)}}{n}\vee\frac
{s\log
(d/s)}{n} \biggr\} \le
\frac{D_2^s \kappa^{-s/4}\sqrt{s\log(d/s)}}{n}\vee\frac{s\log
(d/s)}{n}.
\]
Combining these results, we see that under the conditions $2D_1^*s\log
(d/s)/n\le\kappa$ and
\[
2D_1^*\frac{D_2^s \sqrt{s\log(d/s)}}{n}\le\kappa^{1+{s}/4},
\]
consistent estimation of $J$ is possible. Taking $D^*=2D_1^*(1+D_2^4)$,
we complete the proof of the first claim of the proposition.
To prove the second assertion, we apply Theorem~\ref{thm3}. Since it
holds that $
2\gammaL\ge\gammaL+1\ge\frac{\vartheta}{1+(\sfh(z_{\gammaL
+1})-1)^{-1}}\ge
\frac{\vartheta}{1+(2z_{1})^{-1}}$,
we deduce from Theorem~\ref{thm3} that there are some constants $D_3$
and $D_4$ such that
if
\[
D_3 \biggl\{\frac{\sqrt{N(s,D_4/{\kappa})s\log(d/s)}}{n}\vee \frac
{s\log(d/s)}{n}
\biggr\}\ge\kappa,
\]
then consistent estimation of $J$ is impossible. Since the $s$-dimensional
$L_2$ ball with radius $\sqrt{s\gamma}$ contains the $L_\infty$ ball of
radius $\sqrt{\gamma}$,
$N(s,D_4/{\kappa})\ge(D_5)^s\kappa^{-s/2}$ for some constant $D_5$. By
rearranging different terms, we get
the desired result.
%
\begin{lemma}\label{lemcard0}
For every $\gamma\ge1$ and $d^*\in\NN$, $
N_1(d^*,\gamma)\le0.3 (9\pi e \gamma)^{d^*/2}$.
\end{lemma}
\begin{pf}
One readily checks that if $\|\bk\|_2^2\le d^*\gamma$, then the
hypercube centered at $\bk$ with side of length $1$ is included in the
ball centered at the origin and having
radius $\sqrt{d^*\gamma}+0.5\sqrt{d^*}$. Therefore, $N_1(d^*,\gamma
)\le
(\sqrt{d^*\gamma}+0.5{ \sqrt{d^*}})^{d^*}\operatorname{Vol}[B_{d^*}(0;1)]$,
where $\operatorname{Vol}[B_{d^*}(0;1)]$ stands for the volume of the
unit ball in $\RR^{d^*}$. Using the well-known formula
for the latter and the Stirling approximation, for every $d^*\ge1$, we
get $
\operatorname{Vol}[B_{d^*}(0;1)]= \frac{2\pi^{d^*/2}}{d^*\Gamma(d^*/2)}
\le0.4\frac{(4\pi e/d^*)^{d^*/2}}{\sqrt{2d^*}}$.
This implies that $
N_1(d^*,\gamma)
\le0.4  (\frac{9\gamma d^*}{4} )^{d^*/2}\frac{(4\pi
e/d^*)^{d^*/2}}{\sqrt{2d^*}}
\le\break 0.3 (9\pi e \gamma)^{d^*/2}$ and the result follows.
\end{pf}
\end{appendix}

\section*{Acknowledgments}
The authors would like to thank the reviewers for very useful remarks.

\begin{supplement}
\stitle{Proofs of some results}
\slink[doi]{10.1214/12-AOS1046SUPP} 
\sdatatype{.pdf}
\sfilename{aos1046\_supp.pdf}
\sdescription{The supplementary material provides the proof of
Theorem~\ref{thm1},
Proposition~\ref{prop4},
Lemma~\ref{lem2} and Corollary~\ref{cor2}, as well as those of some
technical lemmas.}
\end{supplement}


%
%

%

\printaddresses

\end{document}